%% file: m3as-review_arxiv.tex
\documentclass[final]{ws-m3as}

\usepackage{epic,eepic}
\usepackage{graphicx,color}
\usepackage[color,notref,notcite]{showkeys}
\usepackage{graphicx}
\usepackage{url}
\usepackage{amssymb,amsmath}
\usepackage{colortbl}
\usepackage{multicol}
\usepackage{multirow}
\usepackage{bm}
\usepackage{citesort}

\definecolor{labelkey}{rgb}{0.6,0,1}

\newcounter{cst}
\catcode`\@=11
\def\ctel#1{C_{\refstepcounter{cst}\@bsphack
\protected@write\@auxout{}%
           {\string\newlabel{#1}{{\thecst}{\thepage}}}\thecst}}
\catcode`\@=12

\newcounter{cexp}
\catcode`\@=11
\def\terml#1{T_{\refstepcounter{cexp}\@bsphack
\protected@write\@auxout{}%
           {\string\newlabel{#1}{{\thecexp}{\thepage}}}\thecexp}}
\catcode`\@=12

\newcommand{\mathbi}[1]{{\boldsymbol #1}}

\newcommand{\ba}{\begin{array}{llll}   }
\newcommand{\bac}{\begin{array}{c}}
\newcommand{\bari}{\begin{array}{r}}
\newcommand{\ea}{\end{array}}

\newcommand{\ban}{\begin{array}{llll}}
\newcommand{\ean}{\end{array}}

\newcommand{\be}{\begin{equation}}
\newcommand{\ee}{\end{equation}}

\newcommand{\beqsys }{\beqtab \left \{ \begin{array}{l}}
\newcommand{\eeqsys }{\end{array} \right . \eeqtab }

\newcommand{\benum}{\begin{enumerate}}
\newcommand{\eenum}{\end{enumerate}}

\newcommand{\beqtab}{\begin{eqnarray}} 
\newcommand{\eeqtab}{\end{eqnarray}}

\newcommand{\dsp}{\displaystyle}

\renewcommand{\d}{{\rm d}}

\newcommand{\disc}{{\cal D}}
\newcommand{\dr}{\partial}

\renewcommand{\div}{{\rm div}}

\renewcommand{\O}{\Omega}
\renewcommand{\phi}{\varphi}

\newcommand{\R}{\mathbb R}
\def\divmfd{\mathcal{DIV}^h}

\newcommand{\bu}{{\overline{u}}}
\def\gradhmm{\mathbf{v}}
\def\n{\mathbf{n}}

\newcommand{\bx}{\overline{\mathbi{x}}}
\def\bxs{\bx_\sigma}
\newcommand{\xcv}{\mathbi{x}}
\def\vect#1{\overrightarrow{#1}}

\def\dist{{\rm d}}

\def\udir{\bu_b}
\def\V{\mathbi{v}}

\begin{document}

\title{FINITE VOLUME SCHEMES FOR DIFFUSION EQUATIONS: INTRODUCTION TO AND REVIEW
OF MODERN METHODS\footnote{Preprint of an article published in
Math. Models Methods Appl. Sci. (M3AS) 24 (2014), no. 8, 1575-1619 (special issue on Recent Techniques for PDE Discretizations on Polyhedral Meshes). DOI:10.1142/S0218202514400041 \copyright{} World Scientific Publishing Company
http://www.worldscientific.com/worldscinet/m3as }}

\author{JEROME DRONIOU}

\address{School of Mathematical Sciences,
Monash University\\
Victoria 3800, Australia.\\
jerome.droniou@monash.edu}

\maketitle

\begin{abstract} We present Finite Volume methods for diffusion equations on generic meshes,
that received important coverage in the last decade or so.
After introducing the main ideas and construction principles of the methods,
we review some literature results, focusing on two important properties
of schemes (discrete versions of well-known properties of the continuous
equation): coercivity and minimum-maximum principles. Coercivity ensures the
stability of the method as well as its convergence under assumptions
compatible with real-world applications, whereas minimum-maximum principles
are crucial in case of strong anisotropy to obtain physically meaningful
approximate solutions.
\end{abstract}



\section{Introduction}\label{sec:intro}

Diffusion processes are ubiquitous in physics of flows, such as heat propagation
or flows in porous media encountered in reservoir engineering. A simple
form of diffusion equation is
\be\label{base}
\ba -\div(\Lambda(x)\nabla\bu(x))=f(x)\,,&\quad x\in\Omega,\\
\bu(x)=\udir\,,&\quad x\in\partial\Omega,
\ea
\ee
where $\Omega$ is the domain of study, $f$ describes the volumic sources or sinks,
$\Lambda$ encodes the diffusion properties of the medium,
$\udir$ is the fixed boundary condition and
$\bu$ is the unknown of interest (pressure, saturation, etc.).
Although very simplified with respect to real-world models, Equation \eqref{base} already
contains some of the main issues that have to be dealt with when designing
and analysing numerical methods for diffusion processes.
The assumptions on the data are:
\begin{eqnarray}
\label{hyp-omega}
&&\Omega\mbox{ is a bounded connected polygonal open subset of $\R^d$, $d\ge 1$,}\\
\label{hyp-fudir}
&&f\in L^2(\Omega)\,,\quad \udir\in H^{1/2}(\O)\,,\\
\label{hyp-lambda}
&&
\ba
\Lambda:\Omega\to \R^{d\times d}\mbox{ is symmetric-valued, essentially bounded and coercive}\\
\mbox{(i.e. $\exists \lambda_-,\lambda_+>0$ such that, for a.e. $x\in\Omega$ and all
$\xi\in \R^d$,}\\
\lambda_- |\xi|^2\le \Lambda(x)\xi\cdot\xi\le \lambda_+|\xi|^2)
\ea
\end{eqnarray}
($\cdot$ and $|\cdot|$ are the Euclidean dot product and norm on $\R^d$).
No other regularity properties are assumed on $\Lambda$, $f$ or $\udir$, and
the proper mathematical formulation of \eqref{base} is therefore,
denoting by $\gamma:H^1(\Omega)\mapsto H^{1/2}(\partial\Omega)$ the trace operator:
\be\label{basew}
\ba
\bu\in \{v\in H^1(\Omega)\;:\;\gamma(v)=\udir\},\\
\dsp\forall \varphi\in H^1_0(\Omega)\,,\quad \int_\Omega \Lambda(x)\nabla\bu(x)\cdot
\nabla\varphi(x)\d x = \int_\Omega f(x)\varphi(x)\d x.
\ea
\ee
Amongst the numerous families of numerical methods for diffusion equations
(Finite Difference, Finite Element, Discontinuous Galerkin...),
Finite Volume (FV) schemes are methods of choice for a number of engineering
applications in which the conservation of various extensive
quantities is important. Local conservativity
of the fluxes is in particular essential to handle the hyperbolicity and
strong coupling which occur in models of miscible
or immiscible flows in porous media.

The purpose of this work is to present a few
modern FV methods for \eqref{base} and to review some of the
mathematical results established for these methods. Although FV
methods can be applied on a number of fluid models, our discussion
will be made with models of porous media flows in mind.
In this case, \eqref{base} corresponds to a steady single-phase
single-component Darcy problem with no gravitational effects,
$\bu$ is the pressure and $\Lambda$ is the permeability field\cite{DIP13-2}.

The paper is organised as follows. In the rest of this section,
we detail the basics behind the construction of FV methods
and we point out two important properties of Equation \eqref{base}
(coercivity and minimum-maximum principle) which are also desirable for discretisations thereof.
Coercivity, in particular, is at the core of techniques which allows
one to carry out convergence proofs without assuming 
non-physical regularities on the data or the solution.
Sec. \ref{sec:FV2} presents the most classical FV method for \eqref{base},
based on a 2-point flux approximation, and highlights its coercivity
and minimum-maximum principle properties as well as its main flaw: it is hardly
applicable on meshes encountered in practical applications.
Secs. \ref{sec:MPFA}, \ref{sec:HMM} and \ref{sec:DDFV} then present
three families of FV schemes applicable on generic meshes: Multi Point
Flux Approximation methods (O-, L- and G-methods), Hybrid Mimetic Mixed methods
(including Hybrid Finite Volume methods, Mimetic Finite Difference schemes and
Mixed Finite Volume methods) and Discrete Duality Finite Volume methods.
In each of these sections, we first present the construction of the method,
focusing on its principles rather than on the details of the
computations, and we then review the literature results on their
coercivity (and convergence) and minimum-maximum principle properties. These sections
are also completed by short conclusions summarising the strengths and weaknesses
of each method. In Sec. \ref{sec:LMP}, we consider some FV
schemes specifically designed to satisfy minimum-maximum principles
on any mesh. Sec. \ref{sec:concl} concludes the paper.

\subsection{What is a Finite Volume scheme?}\label{sec:whatis}

Good question... not easy to answer given the number of methods
presented in the literature as ``Finite Volume'' schemes.
Nevertheless, some basic ideas remain which should be shared by any
method called ``Finite Volume''.

The physical principle that leads to \eqref{base} is the balance of
some extensive quantity $Q$ (heat, component mass, etc.):
given a domain $\omega$, the variation
of $Q$ inside $\omega$ comes from the creation of $Q$ in $\omega$
and the transfer of $Q$ through $\partial\omega$. In a stationary
context, there is no variation of $Q$ and the volumic creation
inside $\omega$ must therefore balance out the quantity of $Q$ which
leaves $\omega$ through $\partial\omega$.
Under modelling assumptions, the creation of $Q$ inside $\omega$
has a volumetric density function $f$ and the flow of $Q$ outside $\omega$
has a surfacic density $-\Lambda(x)\nabla\bu(x)\cdot\n_{\omega}(x)$
(Darcy's or Fourier's law), where $\n_\omega$ is the outer unit normal to $\partial\omega$
and $\Lambda(x)$ is a symmetric positive definite matrix --- heat conductivity matrix in the case
of the heat equation, permeability matrix in reservoir engineering.
The mass balance of $Q$ then reads
\be\label{equation-cons2}
\int_{\partial\omega} -\Lambda(x)\nabla\bu(x)\cdot\mathbf{n}_\omega(x)\d S(x)=\int_\omega f(x)\d x.
\ee
Using Stokes' formula on the left-hand side,
taking $\omega$ a ball around $x\in\O$, dividing by the measure of $\omega$
 and letting its radius tend to $0$ leads to \eqref{base}. This is the ``infinitesimal'' control volume technique to
derive the diffusion equation.

If, on the other hand, we consider a ``finite'' control volume approach in which
$\omega=K$ is a (small but not infinitesimal) polygonal
open set, then \eqref{equation-cons2} becomes
\be\label{balance-flux}
\sum_{\sigma\mbox{ \scriptsize edge of }K} \overline{F}_{K,\sigma}=\int_K f(x)\d x
\ee
where $\overline{F}_{K,\sigma}=\int_\sigma -\Lambda(x)\nabla \bu(x)\cdot \mathbf{n}_K(x)\d S(x)$
is the flux of $\bu$ through $\sigma$.
It can also be noticed that, if $\sigma$ is an edge
between two polygons $K$ and $L$, then
\be\label{cons-flux}
\overline{F}_{K,\sigma}+\overline{F}_{L,\sigma}=0.
\ee

\begin{remark} Another way to get \eqref{balance-flux} is
to integrate \eqref{base} on $K$. This is how FV
methods are usually presented in textbooks, but it is
important to realise that \eqref{balance-flux} directly comes from
physical principles (without even writing \eqref{base}).
This explains why FV methods are particularly
attractive in many engineering contexts.
\end{remark}

The balance \eqref{balance-flux} and conservativity \eqref{cons-flux} of
the fluxes are the two main elements on which FV methods are built.
Let $(\mathcal M,\mathcal E,\mathcal P)$ be a mesh
of $\Omega$ as given by Definition \ref{def:mesh} below.
All FV methods we consider here have at least cell unknowns $(u_K)_{K\in\mathcal M}$,
that play the role of approximate values of $(\bu(\xcv_K))_{K\in\mathcal M}$. 
Such cell unknowns are often desirable in applications, for
coupling issues and because the medium properties (permeability, etc.) are
usually constant in each cell. Some FV methods also use
additional unknowns, e.g. approximate values of $\bu$ on the edges.
The principle of FV schemes is to compute, using all these unknowns,
consistent approximations $F_{K,\sigma}$ of $\overline{F}_{K,\sigma}$
and to write discrete versions of \eqref{balance-flux} and \eqref{cons-flux}:
\be\label{bf}
\mbox{for any }K\in \mathcal M\,:\;\sum_{\sigma\in\mathcal E_{K}} F_{K,\sigma}=\int_K f(x)\d x,
\ee
\be\label{cf}
\mbox{for any edge $\sigma$ between two distinct $K,L\in\mathcal M$}\,:\;
F_{K,\sigma}+F_{L,\sigma}=0.
\ee

\begin{definition}[Mesh]\label{def:mesh} A mesh of $\Omega$ is $(\mathcal M,\mathcal E,\mathcal P)$
where:
\begin{itemize}
\item $\mathcal M$ is a finite family of non-empty open disjoint polygons
(the ``control volumes'' or ``cells'') 
such that $\overline{\Omega}=\cup_{K\in\mathcal M}\overline{K}$,
\item $\mathcal E$ is a finite family of
non-empty disjoint planar subsets of $\Omega$ (the ``edges'') with positive
$(d-1)$-dimensional measure. We assume that
for each control volume $K$ there exists $\mathcal E_K\subset \mathcal E$
such that $\partial K=\cup_{\sigma\in\mathcal E_K}\overline{\sigma}$.
We also assume that each edge $\sigma\in\mathcal E$ belongs to exactly
one or two sets $(\mathcal E_K)_{K\in\mathcal M}$. 
\item $\mathcal P$ is a family of points $(\xcv_K)_{K\in\mathcal M}$ such that,
for each $K$, $\xcv_K\in {K}$.
\end{itemize}
We denote by $|K|$ the $d$-dimensional measure of $K\in\mathcal M$,
by $|\sigma|$ the $(d-1)$-dimensional measure of $\sigma\in\mathcal E$
and by $\mathbf{n}_{K,\sigma}$ the unit normal to $\sigma\in \mathcal E_K$
outward $K$. We also partition $\mathcal E$ into the interior edges $\mathcal E_{\rm int}$
(those included in $\Omega$) and the exterior edges $\mathcal E_{\rm ext}$
(those included in $\partial\Omega$).
The size of the mesh is $h_\mathcal M=\max_{K\in\mathcal M}{\rm diam}(K)$.
We also take $\Lambda_K$ a value of $\Lambda$ in $K$ (e.g. $\frac{1}{|K|}\int_K \Lambda$
or $\Lambda(\xcv_K)$ -- in reservoir applications, $\Lambda$ is constant in each
cell $K$).
\end{definition}

\begin{figure}[h!]
\begin{center}
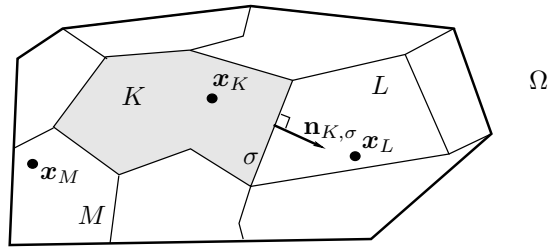
\caption{\label{fig:mesh}A mesh of $\Omega$.}
\end{center}
\end{figure}

\begin{remark} Although we use a 2D vocabulary (polygon, edges...), most of
what we present here is valid in any space dimension.
\end{remark}

\subsection{Convergence analysis and coercivity}\label{sec:analysis}

In reservoir applications, the data (and thus the solution) are not smooth. It is
for example natural for the permeability $\Lambda$ to be discontinuous from
one geological layer to another. Convergence analysis of numerical methods
for such problems should
take into account these practical constraints and should therefore not rely
on non-physical regularity assumptions on the data or solution.
Being able to carry out a convergence analysis under very weak regularity assumptions
on the data or the
solution is also essential for more complex models
(Navier-Stokes equations, multi-phase flows, etc.).

Assuming to simplify that $\udir=0$ (in which case $\bu\in H^1_0(\O)$),
an efficient path to prove the convergence of FV methods for \eqref{base} is to follow
these steps:
\begin{itemize}
\item[\textsf{(C1)}] Establish \emph{a priori} energy estimates on the solution to the scheme,
in a mesh- and scheme-dependent discrete norm which mimics the $H^1_0$ norm,
\item[\textsf{(C2)}] Prove a discrete Rellich compactness result, i.e. that,
as the mesh size tends to $0$, sequences of approximate solutions bounded in these
discrete norms have subsequences which converge(\footnote{In a sense
depending on the method, but which includes at least some form
of strong convergence in $L^2(\Omega)$ and often some form of
weak convergence of discrete gradients.}) to some function $\bu\in H^1_0(\Omega)$,
\item[\textsf{(C3)}] Prove that any such limit $\bu$ of approximate solutions
satisfies \eqref{basew}.
\end{itemize}
Because the solution to \eqref{basew} is unique, Steps \textsf{(C1)}---\textsf{(C3)}
show the convergence of the scheme in the sense that the whole sequence of approximate
solutions converges to the solution of \eqref{basew}. Moreover, for linear schemes,
Step \textsf{(C1)} ensures the existence and uniqueness of a solution to the scheme. 

Following this path does not require any regularity property on $\Lambda$, $f$ or $\bu$
besides those in \eqref{hyp-omega}---\eqref{hyp-lambda} and \eqref{basew}.
Ensuring that \emph{a priori} energy estimates can be obtained
in a proper ``discrete $H^1_0$ norm'' however requires some
assumptions on the scheme.
Consider the continuous equation \eqref{base},
multiply it by $\bu$ and integrate by parts (or, equivalently, take $\varphi=
\bu$ in \eqref{basew}). Then
\begin{equation}\label{ipp}
\lambda_- |\bu|_{H^1_0}^2\le \int_\Omega \Lambda(x)\nabla\bu(x)\cdot\nabla\bu(x)\d x
=\int_\Omega f(x)\bu(x)\d x\le ||f||_{L^2}||\bu||_{L^2}
\end{equation}
and the Poincar\'e inequality $||\bu||_{L^2(\Omega)}\le {\rm diam}(\Omega)|\bu|_{H^1_0(\Omega)}$
gives estimate on $|\bu|_{H^1_0(\Omega)}:=||\,|\nabla\bu|\,||_{L^2(\Omega)}$.
The key element here is the coercivity of $\Lambda$
(which is equivalent to the coercivity of the bilinear form in \eqref{basew}).
Discrete $H^1_0$ estimates on the solution to a FV scheme are usually obtained
by mimicking this process at the discrete level: multiply the scheme by the unknown,
perform discrete integration by parts (or, equivalently, take the unknown
as test function in a variational formulation of the scheme)
and conclude by establishing a discrete Poincar\'e inequality
(see e.g. Sec. \ref{sec:FV2-coer}). This process does not work
for all schemes but, when it does, it shows how to find the
discrete $H^1_0$ norm associated with the scheme and mesh.
For schemes using only cell unknowns, for example, multiplying \eqref{bf} by
$u_K$, summing on the cells and using \eqref{cf}, we see that a
discrete version of \eqref{ipp} can be obtained if there exists
a discrete $H^1_0$ norm $||.||_{1,{\rm disc}}$ satisfying the
Poincar\'e's inequality
\[
\left(\sum_{K\in\mathcal M} |K|u_K^2\right)^{1/2}\le C||u||_{1,{\rm disc}}
\]
and the estimate
\begin{equation}\label{coer-bd}
||u||_{1,{\rm disc}}^2\le C\sum_{\sigma\in\mathcal E} F_{K,\sigma}(u_K-u_L),
\end{equation}
for some $C$ not depending on $u$ or the mesh
(in the previous sum, $K,L$ are the cells on each side of $\sigma\in\mathcal E_{\rm int}$
and $u_L=0$ if $\sigma\in\mathcal E_{\rm ext}\cap\mathcal E_K$).

Obtaining such discrete $H^1_0$ estimates is not only the first step in
proving the convergence of the scheme, but it is also crucial to ensure
its numerical stability. Schemes for which such energy estimates can be established
are called \emph{coercive}. If a linear scheme is coercive and has a symmetric
matrix, then it has a symmetric positive definite matrix and
very efficient algorithms (Cholesky decomposition, conjugate gradient, etc.) can be used to compute
its solution. Note however that the mere symmetry and positive-definiteness
of the matrix are not enough to ensure the coercivity of the scheme, as 
this positive-definiteness must be uniform with respect to
the mesh and must hold for a discrete $H^1_0$ norm satisfying
\textsf{(C2)}.

\begin{remark}[Consistency of Finite Volume methods]
In FV methods, the numerical fluxes $F_{K,\sigma}$ are consistent
approximations of the exact fluxes $\overline{F}_{K,\sigma}$:
if $F^\star_{K,\sigma}$ are the numerical fluxes computed by replacing
the unknowns by the exact values of $\bu$ and if all
data are smooth, then 
\be\label{flux-cons}
F^\star_{K,\sigma}=\overline{F}_{K,\sigma}+\mathcal O(|\sigma|{\rm diam}(K))
\ee
(note that $\overline{F}_{K,\sigma}=\mathcal O(|\sigma|)$).
It is however often said that FV methods do not provide consistent approximations
of the operator $-\div(\Lambda\nabla \bu)$ ``in the Finite Difference sense''
(see Ref. \refcite{EGH00}, Chapter 2).
We can indeed check that, in general,
\be\label{noncons-fd}
\sum_{\sigma\in\mathcal E_K} F^\star_{K,\sigma}=
\int_K -\div(\Lambda\nabla\bu)+\mathcal O(|K|)
\ee
(note that $\int_K -\div(\Lambda\nabla \bu)=\mathcal O(|K|)$).
In fact, as often in mathematical analysis, everything is relative to
topology.
Relation \eqref{noncons-fd} shows a non-consistency in $L^\infty$ or $L^2$
norm, but thanks to the flux consistency \eqref{flux-cons} and the conservativity
of fluxes, we can prove that, for any $\varphi\in H^1_0(\O)$,
\[
\sum_{K\in\mathcal M}\sum_{\sigma\in\mathcal E_K} F^\star_{K,\sigma}\varphi_K
=-\int_\O \div(\Lambda\nabla\bu)(x)\varphi(x)\d x
+\mathcal O(h_{\mathcal M}||(\varphi_K)_{K\in \mathcal M}||_{1,{\rm disc}}),
\]
where $\varphi_K=\frac{1}{|K|}\int_K\varphi(x)\d x$ and
$||\cdot||_{1,{\rm disc}}$ is the discrete $H^1_0$ norm of Sec. \ref{sec:FV2-coer}.
Hence, $\sum_{\sigma\in\mathcal E_K} F^\star_{K,\sigma}$ is a consistent
approximation of $\int_K -\div(\Lambda\nabla \bu)$ in some discrete dual $H^1_0$ norm
and, because of this, establishing discrete $H^1_0$ estimates on approximate
solutions is also crucial to pass to the limit in Step \textsf{(C3)}.
\end{remark}

\begin{remark}[Linearly exact scheme]\label{rem-linex}
The consistency relation \eqref{flux-cons} is strongly related with 
the fact that the scheme is \emph{linearly exact}, meaning that
if the exact solution $\bu$ to \eqref{base} is piecewise linear
on the mesh then its interpolation is the solution to the scheme
(i.e. the scheme exactly reproduces piecewise linear solutions).
In this case, observed numerical orders of convergence(\footnote{Here and everywhere else in this paper, error estimates and orders of convergence are in some form of $L^2$ norm depending on the
scheme.}) are usually 2 for $\bu$ and $1$ for its gradient
(at least for smooth solutions and linear schemes).
\end{remark}

\subsection{Maximum and minimum principles, or monotony}\label{sec:monotony}

A remarkable property of diffusion equations such as \eqref{base} is their
maximum and minimum principles, see Ref. \refcite{HOP27} or Chapter I in Ref. \refcite{MIR70}. In its strong form (also called the local
minimum principle), the minimum
principle states that, should $f$ be non-negative, the solution $\bu$ to \eqref{base}
cannot have a local minimum inside $\Omega$ unless it is constant. This prevents in particular
the solution from presenting oscillating behaviours.
This local minimum principle implies the following weaker (global) form
\begin{equation}\label{minpple}
\mbox{if $f\ge 0$ and $\udir\ge 0$ then $\bu\ge 0$},
\end{equation}
as well as the (global) minimum-maximum principle (obtained
by applying \eqref{minpple} to $\bu-(\inf_{\dr\O}\udir)$ and $(\sup_{\dr \O} \udir)-\bu$):
\begin{equation}\label{minmaxpple}
\mbox{if $f=0$ then $\inf_{\dr\O} \udir \le \bu\le \sup_{\dr\O}\udir$.}
\end{equation}

Assume that $U=((u_i)_{i\in I},(u_z)_{z\in B})$ is a vector gathering
the unknowns $(u_i)_{i\in I}$ of the scheme and the discretised
boundary conditions $(u_z)_{z\in B}$, computed from $\udir$.
If the scheme is written $S(U)=R$, where $R$ is
a vector constructed from $f$, the discrete desirable versions of \eqref{minpple}
and \eqref{minmaxpple} are
\begin{equation}\label{disc-minpple}
\mbox{if $S(U)=R \ge 0$ and $u_z\ge 0$ for all $z\in B$ then $u_i\ge 0$
for all $i\in I$}
\end{equation}
(where $R\ge 0$ means that all components of $R$ are non-negative)
and
\begin{equation}\label{disc-minmaxpple}
\mbox{if $S(U)=0$ then $\inf_{z\in B} u_z\le u_i
\le \sup_{z\in B}u_z$ for all $i\in I$}.
\end{equation}
For linear schemes (i.e. $S$ is a linear function) that are
exact on constant functions (i.e. $S(\mathbf{1})=0$, where
$\mathbf{1}$ is the vector with all components equal to $1$),
the discrete minimum principle \eqref{disc-minpple} implies
the discrete minimum-maximum principle \eqref{disc-minmaxpple}
(if $S(U)=0$, apply \eqref{disc-minpple} to $V=(\max_{z\in B}u_z) \mathbf{1} - U$ 
and $V=U-(\min_{z\in B}u_z)\mathbf{1}$, which both satisfy
$V_b\ge 0$ for all $b\in B$ and $S(V)=0$ by linearity of $S$).
As we shall see in Sec. \ref{sec:LMP}, non-linear schemes may satisfy
\eqref{disc-minpple} without satisfying
\eqref{disc-minmaxpple}.

The usual way in the literature to prove that a linear
scheme satisfies \eqref{disc-minpple} is to show that its
matrix $A=(a_{ij})_{ij}$ is diagonally dominant by columns
(i.e. $a_{ii}>0$ for all $i$, $a_{ij}\le 0$ for all $i\not= j$ and
$a_{kk}\ge \sum_{i\not=k}|a_{ik}|$ for all $k$ with strict inequality
for at least one $k$) and has a connected graph.
Under these assumptions, it is easy to see that
$A$ is invertible and that $A^{-1}$ only has non-negative coefficients
($A$ is thus an $M$-matrix), see Chapter 6 in Ref. \refcite{ABR79}.
Provided that the scheme is written
$S(U)=A(u_i)_{i\in I} - C(u_z)_{z\in B}=R$ where $C$ is a matrix
with non-negative coefficients, we then obtain $(u_i)_{i\in I}=A^{-1}(R+C(u_z)_{z\in B})
\ge 0$ whenever $R\ge 0$ and $u_z\ge 0$ for all $z\in B$.

Satisfying a discrete minimum-maximum principle is particularly important
in complex models such as multi-phase flows in reservoir engineering.
Schemes that do not satisfy this principle may give rise
to spurious oscillations which may lead to gas-oil numerical
instabilities.
Linear schemes for \eqref{base} satisfying \eqref{disc-minpple}
are also called \emph{monotone}, as they preserve the order
of boundary conditions (for non-negative right-hand sides) or of
initial conditions (when applied to transient equations).

\section{TPFA scheme}\label{sec:FV2}

Let us assume that the medium is isotropic, i.e. $\Lambda(x)=\lambda(x){\rm Id}$
for some scalar function $\lambda$. We also assume the following orthogonality conditions
on the mesh:
\be\label{cond-orth}
\ba
\forall \sigma\mbox{ edge between two control volumes $K,L\in\mathcal M$}\,,\;
(\xcv_K\xcv_L)\bot \sigma\,,\\
\forall \sigma\in\mathcal E_{\rm ext}\cap\mathcal E_K\,,\;
\mbox{ the half-line $\xcv_K+[0,\infty)\mathbf{n}_{K,\sigma}$
intersects $\sigma$}.
\ea
\ee
In Fig. \ref{fig:mesh}, for example, this assumption is satisfied by
the edge $\sigma$ between $K$ and $L$ but not by the
edge between $K$ and $M$.
Letting $\{\xcv_\sigma\}=(\xcv_K\xcv_L)\cap \sigma$ (or $\{\xcv_\sigma\}=(\xcv_K+[0,\infty)\mathbf{n}_{K,\sigma})
\cap \sigma$ if $\sigma\in\mathcal E_{\rm ext}$), consistent approximations
of the fluxes for small $h_{\mathcal M}$ are
\begin{eqnarray}
\label{disc-flux-int}
\mbox{if $\sigma\in\mathcal E_K\cap\mathcal E_L$}&:&
F_{K,\sigma}= \lambda_K|\sigma|\frac{u_K-u_\sigma}{\dist(\xcv_K,\xcv_\sigma)}
\mbox{ and } F_{L,\sigma}= \lambda_L|\sigma|\frac{u_L-u_\sigma}{\dist(\xcv_L,\xcv_\sigma)},\\
\label{disc-flux-ext}
\mbox{if $\sigma\in\mathcal E_{\rm ext}\cap \mathcal E_K$}
&:&F_{K,\sigma}=\lambda_K|\sigma|\frac{u_K-u_\sigma}{\dist(\xcv_K,\sigma)},
\end{eqnarray}
where $\dist(a,b)=|a-b|$, $\lambda_K$ is the value of $\lambda$ on $K$
and $u_\sigma$ approximates $\bu(\xcv_\sigma)$.
If $\sigma\in\mathcal E_{\rm ext}$, $u_\sigma$ is fixed by $\udir$(\footnote{Several
choices are possible. If $\udir$ is smooth enough, then one can take
$u_\sigma=\udir(\xcv_\sigma)$.
Otherwise, $u_\sigma$ can be chosen as the average of $\udir$ on $\sigma$.}).
If $\sigma\in\mathcal E_K\cap\mathcal E_L$, the
additional unknown $u_\sigma$ is eliminated by imposing the conservativity
\eqref{cf} of fluxes and we get (see Ref. \refcite{EGH00}, Chapter 3):
\be\label{2pt-fluxes-h}
F_{K,\sigma}=\tau_\sigma (u_K-u_L) \quad\mbox{ with }
\tau_\sigma
=\frac{|\sigma|}{\dist(\xcv_K,\xcv_L)}\frac{\lambda_K \lambda_L \dist(\xcv_K,\xcv_L)}{\lambda_K\dist(\xcv_L,\xcv_\sigma)
+\lambda_L  \dist(\xcv_K,\xcv_\sigma)}.
\ee
The balance equation \eqref{bf} of the discrete fluxes \eqref{disc-flux-ext}-\eqref{2pt-fluxes-h}
then gives an FV scheme for \eqref{base} when $\Lambda=\lambda{\rm Id}$,
called the Two Point Flux Approximation Finite Volume scheme (TPFA for short)
since each flux is computed using only the 2 unknowns on each side of the
edge.

\begin{remark} As $\dist(\xcv_K,\xcv_\sigma)+\dist(\xcv_L,\xcv_\sigma)=\dist(\xcv_K,\xcv_L)$,
the transmissibility $\tau_\sigma$ involves an harmonic average
of the values of $\Lambda$ in the cells on each side of $\sigma$. This harmonic
average is well-known, in FV methods, to give a much more accurate solution
than other averages.
\end{remark}

\begin{remark}\label{rem:orth} If $\Lambda$ is an anisotropic full tensor,
the same construction can be made (see Chapter 3 in Ref. \refcite{EGH00}) provided that
the orthogonality condition \eqref{cond-orth} is replaced with \eqref{cond-orth2},
in which $D_{K,\sigma}$ is the straight line going through $\xcv_K$ and orthogonal to $\sigma$
for the scalar product induced by $\Lambda_K^{-1}$:
\be\label{cond-orth2}
\ba
\forall \sigma\mbox{ between two control volumes $K,L\in\mathcal M$}\,,\;
D_{K,\sigma}\cap \sigma=D_{L,\sigma}\cap\sigma\not=\emptyset\,,\\
\forall \sigma\in\mathcal E_{\rm ext}\cap\mathcal E_K\,,\;
D_{K,\sigma}\cap\sigma\not=\emptyset.
\ea
\ee
\end{remark}

\subsection{Coercivity}\label{sec:FV2-coer}

Assume that $\udir=0$ and thus that $u_\sigma=0$ for all $\sigma\in\mathcal E_{\rm ext}$.
Multiplying the balance equation \eqref{bf} by $u_K$,
summing on $K\in\mathcal M$ and gathering by edges (=discrete integration
by parts), we obtain, thanks to \eqref{2pt-fluxes-h},
\be\label{2pt-estimate}
||u||_{1,\disc}^2:=\sum_{\sigma\in\mathcal E_{\rm int}}\tau_\sigma(u_K-u_L)^2+
\sum_{\sigma\in\mathcal E_{\rm ext}}\tau_\sigma u_K^2
=\int_\Omega f(x)u(x)\d x
\ee
where $u$ is the piecewise constant function equal to $u_K$ on $K$
and, in the sums, $K$ and $L$ are the control volumes
on each side of $\sigma\in\mathcal E_{\rm int}$
(we let $\tau_\sigma=\lambda_K\frac{|\sigma|}{\dist(\xcv_K,\sigma)}$
whenever $\sigma\in\mathcal E_{\rm ext}\cap\mathcal E_K$).
The left-hand side of \eqref{2pt-estimate} defines a discrete
$H^1_0$ norm $||u||_{1,{\rm disc}}$ for which one can
establish the discrete Poincar\'e inequality
$||u||_{L^2(\Omega)}\le {\rm diam}(\Omega)||u||_{1,{\rm disc}}$
and a discrete compactness result as in  Step \textsf{(C2)}
of Sec. \ref{sec:analysis}, see Chapter 3 in Ref. \refcite{EGH00}.
The TPFA scheme is thus coercive (with a symmetric
matrix) and its convergence can be proved under the sole
assumptions \eqref{hyp-omega}--\eqref{hyp-lambda}. Of course,
error estimates can also be obtained if the data are more regular\cite{HER95}.

\subsection{Monotony}

Injecting \eqref{disc-flux-ext}-\eqref{2pt-fluxes-h} in the balance equation
\eqref{bf} we obtain, with the same conventions as in \eqref{2pt-estimate},
for all $K\in\mathcal M$,
\be\label{2pt-mpple}
\sum_{\sigma\in\mathcal E_{\rm int}}\tau_\sigma(u_K-u_L)
+\sum_{\sigma\in\mathcal E_{\rm ext}}\tau_\sigma u_K
=\dsp \int_K f(x)\d x
+\sum_{\sigma\in\mathcal E_{\rm ext}}\tau_\sigma u_\sigma.
\ee
{}From this expression we can see that the scheme's function (see Section \ref{sec:monotony})
can be written $S(U)=A(u_K)_{K\in \mathcal M}-C(u_\sigma)_{\sigma\in\mathcal E_{\rm ext}}$,
with $A$ diagonally dominant, symmetric and graph-connected, and all coefficients of $C$
non-negative. Sec. \ref{sec:monotony} then shows that the TPFA scheme is monotone.

\begin{remark}\label{rem:monvf2}
Monotony of the TPFA scheme is in fact easy to prove from \eqref{2pt-mpple}.
If $f\ge 0$, $u_\sigma\ge 0$ for all $\sigma\in\mathcal E_{\rm ext}$
and $u_K=\min_{M\in\mathcal M}u_M<0$ then
the left-hand side of \eqref{2pt-mpple} is a 
non-negative sum of non-positive terms. Hence all terms are equal to $0$
and $u_K=u_L$ for all neighbours $L$ of $K$. The minimal value $u_K$ 
thus propagates to all neighbours and, ultimately, to the whole connected domain.
Using \eqref{2pt-mpple} for one boundary cell
then contradicts the negativity of this minimal value.

In fact, this reasoning applied to $A^T$ gives a proof that the diagonal dominance
by columns of $A$ and its graph connectedness entail the non-negativity of all coefficients
of $A^{-1}$. It also shows that schemes with such matrices
satisfy in fact a discrete version of the strong minimum principle:
if $f\ge 0$, the solution to the scheme cannot have any interior minimum
unless it is constant.
\end{remark}

\subsection{The perfect scheme?}

The TPFA scheme is a cell-centred scheme (only involving cell unknowns), very
cheap to implement and with a small stencil:
5 on 2D quadrilateral meshes and 7 on 3D hexahedral meshes. Its matrix is
therefore very sparse and its solution easy to compute.
For these reasons, it has been adopted in many engineering software,
but it is not the perfect scheme...

Meshes available in field applications
may be quite distorted and may have cells presenting various complex geometries,
especially in basin simulation where alignment with
geological layers and erosion may lead to hexahedra with
collapsed faces. The orthogonality properties
\eqref{cond-orth} or \eqref{cond-orth2} are impossible to satisfy
on these meshes and, should they fail for too many edges, the solution given
by the TPFA scheme will be totally incorrect\cite{FAI92,AAV02,EIG05}.
Other FV methods therefore had to be designed, providing consistent fluxes
for general meshes and tensors.

\section{MPFA methods}\label{sec:MPFA}

Consistent approximations of the fluxes $\overline{F}_{K,\sigma}$
on general meshes require the usage of more approximate values of $\bu$
(in cells, on edges or at vertices) than the two at $\xcv_K$
and $\xcv_L$ on each side of $\sigma$.
One easy way to get such values is to interpolate them from cell unknowns.
This is the path chosen in Ref. \refcite{FAI92} which introduces,
for each edge, additional cell values located at points satisfying the
orthogonality condition \eqref{cond-orth2} for the considered edge,
and then compute these values by convex combinations of existing cell unknowns.
However, this scheme's construction and stability can only be
ensured for grids not too distorted and tensors not too anisotropic.

Another idea is not to try and get back the orthogonality condition
\eqref{cond-orth2}, but to use the additional values to compute
approximate gradients, which in turn give approximate fluxes $F_{K,\sigma}$.
However, the computation of the
additional values must be done in a clever way, especially when $\Lambda$
is discontinuous, to ensure that the flux conservativity \eqref{cf} is satisfied.

The Multi-Point Flux Approximation (MPFA) schemes are based on such a
construction. Introduced in the mid- to late 90's\cite{AAV96,AAV98-I,AAV98-II,EDW98,EDW94},
these methods assume that the solution is piecewise linear in some sub-cells around each
vertex, introduce additional edge unknowns and express the linear variation
of the solution to compute gradients and thus fluxes in these sub-cells.
The edge unknowns are then eliminated (interpolated using cell unknowns) by writing
\emph{continuity equations} for the solution and
\emph{conservativity equations} for its fluxes. The final numerical fluxes
are consistent, conservative and expressed only in terms of cell unknowns.

\subsection{O-method}

Several MPFA methods have been devised over the years and their main variation is
on the choice of the local continuity and conservativity equations.
Amongst those methods, the O-method
(presented in Refs. \refcite{AAV02,AAV98-I} for particular polygonal
meshes) has received one of the largest coverage in literature
on MPFA methods.

\begin{figure}[h!]
\begin{center}
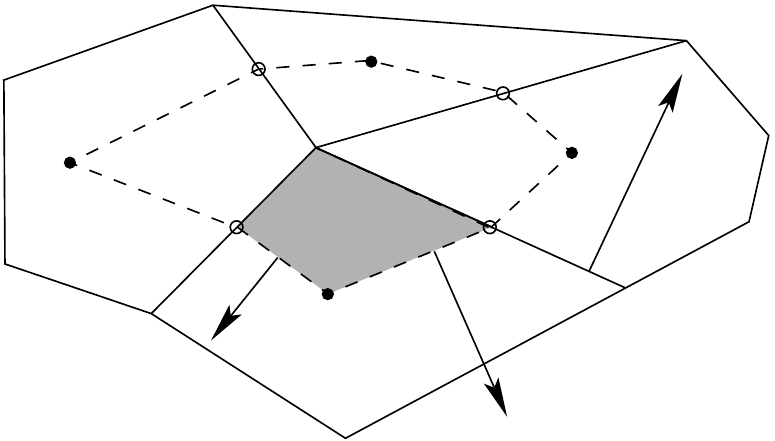
\caption{\label{fig:mpfa}Control volumes ($K,L,\ldots$) and
interaction region (enclosed in dotted line) for the MPFA O-method.
$\nu_{K,\tau}=$ normal
vector to $(\xcv_K\bx_\tau)$ with length $\dist(\xcv_K,\bx_\tau)$ ($\tau=\sigma,\sigma'$).
}
\end{center}
\end{figure}
Let us first consider the 2D case. For each edge $\sigma$, we fix a point $\bxs$
on $\sigma$. Several choices are possible\cite{AAV98-I,EDW98} but we only consider
here the case where $\bxs$ is the midpoint of $\sigma$. Then for each
vertex $\V$ of the mesh, an \emph{interaction region} is built by joining the
cell points $\xcv_K$ around $\V$ and the midpoints $\bxs$ of the edges containing $\V$
(see Fig. \ref{fig:mpfa}). This interaction region is made of one sub-cell $K_\V$ per
cell $K$ and the solution $\bu$ is approximated by
a function that is linear inside each sub-cell around $\V$ (\footnote{This
linear approximation is natural if the mesh size is small enough since,
usually, $\Lambda$ and $f$  are assumed to be constant or smooth in $K$,
so that $\bu$ is expected to be smooth inside $K$.}).

At this stage, \emph{continuity of this piecewise linear approximation is assumed at
each edge midpoint $\bxs$} around $\V$. We can therefore talk about
the value $u_\sigma$ of this function at $\bxs$, and its
constant gradient $\nabla_{K_\V}u$ on $K_\V$ satisfies
\be\label{lin-grad}
\nabla_{K_\V}u\cdot (\xcv_K-\bx_\tau)=u_K-u_\tau\quad \mbox{($\tau=\sigma$ or $\sigma'$)}.
\ee
Assuming that the vectors $\vect{\xcv_K\bxs}$ and $\vect{\xcv_K\bx_{\sigma'}}$ are
linearly independent, these two projections of $\nabla_{K_\V}u$ on these
vectors provide\cite{AAV98-I} the whole gradient $\nabla_{K_\V}u$:
\be\label{eq-grad}
\nabla_{K_\V}u = -\frac{1}{2T}\left((u_\sigma-u_K)\nu_{K,\sigma'}+(u_{\sigma'}-u_K)\nu_{K,\sigma}\right)
\ee
where $T$ is the area of triangle $(\xcv_K\bxs\bx_{\sigma'})$ and
$\nu_{K,\tau}$ ($\tau=\sigma$ or $\sigma'$) is the normal to $(\xcv_K\bx_\tau)$,
pointing outward this triangle and having length $\dist(\xcv_K,\bx_\tau)$.

Sub-fluxes across the half-edges $[\V\bx_\tau]$ around $\V$ are then computed
using these gradients, and therefore depend
on the cell unknowns $u_K,u_L,\ldots$ and the edge unknowns $u_\sigma,u_{\sigma'},\ldots$
around $\V$. For example, the sub-flux from $K$ through
$[\V\bx_\tau]$ is
\be\label{eq-flux}
F_{K,\tau,\V}=-\dist(\V,\bx_\tau)\Lambda_K\nabla_{K_\V}u\cdot\mathbf{n}_{K,\tau}
\quad  (\tau=\sigma\mbox{ or }\sigma'). 
\ee
The next step is to eliminate the edge unknowns involved in these sub-fluxes.
This is done by imposing the conservativity of the fluxes around $\V$:
\be\label{elim-O}
\ba
\mbox{For any edge $\tau$ containing $\V$, if $R,S$ are the
cells on each side of $\tau$,}\\
\dsp F_{R,\tau,\V}+F_{S,\tau,\V}=0.
\ea
\ee
Note that if $\tau$ is an edge on $\partial\Omega$, $u_\tau$
is not eliminated but fixed by the value of $\udir$
(Neumann boundary conditions are also easily handled,
either by imposing the value of $F_{K,\tau,\V}$ whenever $\tau$ is a boundary edge
or by using -- which is equivalent -- ghosts cells outside $\Omega$\cite{AAV98-I,AAV02}).

{}From the construction \eqref{eq-grad}-\eqref{eq-flux} of the sub-fluxes,
\eqref{elim-O} gives a linear square system on the edge unknowns
$u_\sigma,u_{\sigma'},\ldots$ around $\V$ which is, in general, invertible and gives an
expression of these edge unknowns in terms of the cell unknowns $u_K,u_L,\ldots$ around $\V$.
Plugged into \eqref{eq-grad}-\eqref{eq-flux}, these expressions of the edge unknowns
give formulas for the sub-flux $F_{K,\sigma,\V}$ using only the cell unknowns
$u_K,u_L,\ldots$ around $\V$. The same procedure performed
from the other vertex $\V'$ of $\sigma$ gives a second sub-flux $F_{K,\sigma,\V'}$.
The global flux through $\sigma$,  that is
$F_{K,\sigma}=F_{K,\sigma,\V}+F_{K,\sigma,\V'}$,
is therefore a function of all the unknowns $u_K,u_L,\ldots$
in all the cells around $\V$ and $\V'$.
By construction, $(F_{K,\sigma})_{K\in\mathcal M\,,\;\sigma\in\mathcal E_K}$
naturally satisfy the conservativity
equation \eqref{cf} and the O-scheme is thus obtained by only imposing
the balance equation \eqref{bf}.

\begin{remark}[Two edge unknowns per edge] The elimination of the edge unknowns
is performed locally around each vertex $\V$ and the continuity at the edge midpoints is
only enforced when eliminating the edge unknowns around $\V$.
The edge unknown $u_\sigma$ at $\bxs$ when viewed from vertex $\V$ therefore may
be different from the edge
unknown at $\bxs$ viewed from the other vertex $\V'$ of $\sigma$.
This may look strange, as there is no particular reason for
$\bu$ to have different values at $\bxs$, but this comes from the
construction of the MPFA method which cannot assume that the
linear variations of $\bu$ in $K_\V$ and in $K_{\V'}$ have the same
value at $\bxs$ (otherwise, some flux conservativity equations could
not be satisfied).
\end{remark}

The generalisation of this construction to 3D polyhedral cells is pretty straightforward\cite{AAV06}
if we assume that
\be\label{assum-3D}
\mbox{for each cell $K$ and each vertex $\V$ of $K$,
exactly 3 faces of $K$ meet at $\V$.}
\ee
In this case, the sub-cell
$K_\V$ is the hexahedron obtained by joining $\V$, $\xcv_K$, the midpoints of edges
of $K$ having $\V$ as vertex and the three centres of gravity $\bxs$,
$\bx_{\sigma'}$ and $\bx_{\sigma''}$ of the faces of $K$
meeting at $\V$. Three temporary unknowns $u_\sigma$, $u_{\sigma'}$ and $u_{\sigma''}$
are introduced at the centres of gravity of the faces and,
assuming that the vectors $\vect{\xcv_K\bxs}$, $\vect{\xcv_K\bx_{\sigma'}}$ and
$\vect{\xcv_K\bx_{\sigma''}}$ are linearly independent, the three equations
\eqref{lin-grad} for $\tau=\sigma$, $\sigma'$ and $\sigma''$ can be solved
for $\nabla_{K_\V}u$, which is thus computed in terms of $u_K,u_\sigma,u_{\sigma'}$
and $u_{\sigma''}$.
The rest of the construction follows as in 2D, the edge
unknowns being eliminated thanks to the sub-fluxes conservativity.

\begin{remark}
This procedure even allows for non-planar faces (which often occurs
in hexahedral meshes in 3D, as the four vertices of a given face may
not be on the same plane), provided that the vectors $\mathbf{n}_{K,\sigma}$
are defined as the mean value on $\sigma$ of the pointwise normal vector to the face\cite{AAV02,AAV06}.
\end{remark}

Construction of an MPFA O-method on 3D meshes is much less obvious
when \eqref{assum-3D} does not hold. In this case,
for some vertices $\V$ the system \eqref{lin-grad} has 4 or more equations and,
since (in general) the gradient $\nabla_{K_\V}u$ is entirely determined by $u_K$ and
only 3 face unknowns, the other face unknowns will be fixed
by those 3 face unknowns. No degrees of freedom then remain to impose
the conservativity of the corresponding sub-fluxes.
Ref. \refcite{AGE10} however introduces a scheme on general polygonal
or polyhedral meshes (without assuming \eqref{assum-3D}),
which coincides with the MPFA O-method in 2D and in 3D when
\eqref{assum-3D} holds. This reference also presents a new formulation of the O-method,
based on a discrete form of the variational formulation \eqref{basew}
rather than on a flux balance \eqref{bf}.

\begin{remark}
Explicit formulas for the fluxes in terms of the cell unknowns can
be obtained\cite{AAV02} in the case of parallelogram or parallelepiped meshes
and $\Lambda$ constant. In other cases, System \eqref{elim-O} has to be numerically solved.
\end{remark}

\begin{remark} 
For non-conforming meshes such as the ones appearing in reservoirs with
faults, this MPFA O-method leads to unacceptable fluxes and must
therefore be modified\cite{AAV01},
by introducing two linear approximations of $\bu$ in some sub-cells $K_\V$.
\end{remark}

\subsection{L- and G-methods}\label{sec:LG}

As already mentioned, many choices are available to compute consistent
conservative fluxes from piecewise linear approximations of $\bu$ around each vertex.
Another well-studied MPFA method is the L-method,
introduced in Ref. \refcite{AAV08} for quadrilateral meshes.
The major difference of the L-method with respect to the
O-method are: (i) no edge unknowns need to be introduced as
the gradient themselves are the additional unknowns to eliminate,
(ii) the continuity and sub-flux conservativity equations are written
only on 2 edges, (iii) the continuity of the piecewise linear approximation
is imposed on whole edges (not only at edge midpoints),
and (iv) the gradients and piecewise linear approximation constructed on sub-cells $K_\V$, $L_\V$,
$\ldots$, depend on the edge $\sigma$ through which we want to compute the flux
and are thus not common to all sub-fluxes around $\V$.

Still using the notations in Fig. \ref{fig:mpfa},
let us consider the sub-flux $F_{K,\sigma,\V}$ and let us introduce
$\nabla^\sigma_{M_\V}u$, $\nabla^\sigma_{K_\V}u$ and $\nabla^\sigma_{L_\V}u$,
the three constant gradients of a piecewise linear approximation of $\bu$
on $M_\V\cup K_\V\cup L_\V$. As mentioned above, these gradients will only
be used to compute $F_{K,\sigma,\V}$ and other gradients would be used
if we were to compute $F_{K,\sigma',\V}$ for example (ergo the
super-script $\sigma$ in $\nabla^\sigma_{M_\V}u$, $\nabla^\sigma_{K_\V}u$, $\nabla^\sigma_{L_\V}u$).
In the L-method, full continuity is imposed for this approximation:
\be\label{full-cont}
\ba
\dsp \forall x\in [\V\bx_{\sigma'}]\,:\;u_K+(\nabla^\sigma_{K_\V}u)\cdot(x-\xcv_K)=
u_M+(\nabla^\sigma_{M_\V}u)\cdot(x-\xcv_M)\\
\dsp \forall x\in [\V\bxs]\,:\;u_K+(\nabla^\sigma_{K_\V}u)\cdot(x-\xcv_K)
=u_L+(\nabla^\sigma_{L_\V}u)\cdot(x-\xcv_L)
\ea
\ee
These equations can be equivalently written only at $\V$, $\bx_{\sigma'}$
and $\V$, $\bxs$ respectively, and they provide 4 conditions on the 6 degrees of freedom of
the 3 gradients. The sub-flux conservativities give the remaining 2 equations
\be\label{cons-L}
\ba
\dsp \Lambda_K \nabla^\sigma_{K_\V}u\cdot\mathbf{n}_{K,\sigma'}
+\Lambda_M \nabla^\sigma_{M_\V}u\cdot\mathbf{n}_{M,\sigma'}=0\\
\dsp \Lambda_K \nabla^\sigma_{K_\V}u\cdot\mathbf{n}_{K,\sigma}
+\Lambda_L \nabla^\sigma_{L_\V}u\cdot\mathbf{n}_{L,\sigma}=0.
\ea
\ee
System \eqref{full-cont}-\eqref{cons-L} is therefore square
and invertible in general (otherwise, work\-arounds can be
designed\cite{AGE10-2}). The local gradients can then be
expressed in terms of the cell unknowns $u_M$, $u_K$ and $u_L$,
and so does the sub-flux $F_{K,\sigma,\V}=-\dist(\V,\bxs)\Lambda_K \nabla^\sigma_{K_\V}u\cdot
\mathbf{n}_{K,\sigma}$.

\begin{remark}
If $\sigma$ or $\sigma'$ is a boundary edge, then the corresponding
right-hand side in \eqref{full-cont} is fixed by the value of $\udir$ and
the corresponding conservativity equation in \eqref{cons-L} is removed.
System \eqref{full-cont}-\eqref{cons-L} remains square, of size 4
(if only one edge is a boundary edge) or 2 (if both $\sigma$ and $\sigma'$ are boundary edges).
\end{remark}

This is however but one choice that can be made to compute the flux through $\sigma$.
Another natural choice would be to use the edges $\sigma$ and $\sigma''$ instead of $\sigma$
and $\sigma'$ in \eqref{full-cont}-\eqref{cons-L}. This would give another
sub-flux $F_{K,\sigma,\V}$ in terms of $u_K$, $u_L$, $u_N$.
In the L-method, the choice between using $\sigma,\sigma'$ or $\sigma,\sigma''$
is made according to a criterion\cite{AAV08} involving transmissibility signs and ensuring
that each cell unknown $u_M,u_K,u_L$ or $u_K,u_L,u_N$ contributes with the most
physically-relevant sign to the sub-flux through $[\V\bxs]$.
Full formulas can be obtained\cite{AAV08} in the case of homogeneous
media and grids made of parallelograms and, in the case
of moderate skewness of the diffusion tensor and the grid, the chosen
criterion indeed leads to the correct signs.

\begin{remark}The L-method does not suffer from the same issues
(and does not need modification) as the original
MPFA O-method on meshes with faults\cite{AAV08}.\end{remark}
 
A generalisation of the L-method, the G-method, has been proposed in Ref. \refcite{AGE10-2}.
Its principles are the same (full continuity of $\bu$ and conservativity of the
fluxes on some edges), but the above selection criterion is not applied and
the global fluxes through $\sigma$ are built as convex combinations of all possible
sub-fluxes through this edge. These combinations are chosen according to
some local index, designed to improve the coercivity properties of the
scheme.

\begin{remark} Contrary to the O-scheme, construction of
the L- and G-scheme on general 3D polyhedral meshes is straightforward\cite{AGE10-2}.
Indeed, no face unknown is introduced and there is always, whatever the number of faces
that meet at a given vertex, enough degrees of freedom (one local constant gradient
per face which contains the vertex) to impose the local conservativity of sub-fluxes.
\end{remark}

\begin{remark}
The MPFA U-method\cite{AAV98-I} is based on principles a bit similar to the
L-method, computing the flux through $[\V\bxs]$ by mixing
midpoint continuity \eqref{lin-grad} (at $\bxs$) and the full continuity
on $[\V\bx_{\sigma'}]$ and $[\V\bx_{\sigma''}]$ (as in \eqref{full-cont}).
The local gradients also depend on the edge $\sigma$ through which we compute
the sub-flux.
\end{remark}

\subsection{Coercivity and convergence of MPFA methods}

MPFA methods are linearly exact, and therefore consistent in the sense \eqref{flux-cons},
but they are not coercive in general. Using reference elements (or curvilinear coordinates)
such as in Finite Element methods, constructions of symmetric definite positive MPFA O-methods have been
proposed on quadrilateral (hexahedral in 3D) meshes in Refs. \refcite{AAV02,AAV07,EDW08-II}
and on general 2D polygonal meshes in Ref. \refcite{FRI08}.
However, these methods method turn out to be numerically less stable than the
MPFA O-method presented above\cite{AAV07} (constructed in physical space).
Convergence of these reference element-based O-methods even sometimes seems to
be lost in presence of anisotropy or perturbed mesh,
when the O-method constructed in physical space still converges\cite{AAV06,AAV07,KLA06}. A reason for this
loss of convergence, in view of Sec. \ref{sec:analysis},
is probably the following\cite{AAV07}: when constructing the method on a reference
mesh, the coercivity properties of the scheme matrix depends on the
mesh regularity (via the Piola mapping) and may
degenerate for strongly perturbed meshes as the mesh size tends to $0$,
thus preventing from establishing energy estimates in a proper discrete $H^1_0$ norm
for which the compactness result of Step \textsf{(C2)} in Sec. \ref{sec:analysis}
would hold.

It has been proved that the physical O-method is coercive (and
gives a symmetric definite positive matrix) on meshes made of parallelograms (parallelepiped in 3D)
with $(\xcv_K)_{K\in\mathcal M}$ the centres of gravity of the cells\cite{AAV06,AGE10}.
This is also true for meshes made of triangles (tetrahedra in 3D),
provided that the unknown $u_\sigma$ used to construct
the piecewise linear approximation of $\bu$ in $K_\V$ is
not located at $\bxs$ but closer to $\V$ (see Refs. \refcite{AGE10,LPO05}).
Except in those particular instances,
proofs of convergence of MPFA methods are always done by \emph{assuming}
some coercivity property. 

Ref. \refcite{KLA06} compares the MPFA O-method on 2D quadrilateral
meshes to a non-symmetric Mixed Finite Element method (using a particular quadrature
rule) and obtains, under a global coercivity assumption on the system matrix,
$\mathcal O(h_{\mathcal M})$ error estimates for the approximate
solution and fluxes, under the assumptions $\Lambda\in
C^1(\overline{\Omega})$ and $\bu\in H^2(\Omega)$.
In a recent study\cite{KLA12}, the MPFA O-method is compared on 2D or
3D polyhedral meshes satisfying \eqref{assum-3D} to
a non-symmetric Mimetic Finite Difference method (see Sec. \ref{sec:HMM}).
Under local coercivity assumptions, $\mathcal O(h_{\mathcal M}^{\alpha})$
error estimates are obtained when $\Lambda\in C^1(\overline{\Omega})$
and $\bu\in H^{1+\alpha}(\Omega)$ ($\alpha> 1/2$ in 3D).

The regularity assumptions on $\Lambda$ and $\bu$ required to establish
these error estimates are not compatible
with usual field applications (see Sec. \ref{sec:analysis}).
It is  however possible to perform the full convergence analysis of
the MPFA O- and L-method without assuming any non-physical smoothness
on the data, by following the path sketched in Sec. \ref{sec:analysis}.
This is done in Ref. \refcite{AGE10} for the MPFA O-method and in Ref. \refcite{AGE10-2}
for the MPFA L- and G-method. In these references, the convergence of MPFA
methods on generic grids, in 2D or 3D (without assuming \eqref{assum-3D}),
is proved by only assuming \eqref{hyp-omega}---\eqref{hyp-lambda} and some
local coercivity conditions which can be checked in numerical experiments.

The numerical study of the convergence of MPFA methods has also been performed in
a number of articles\cite{EIG05,AAV06,PAL06}. As
expected, the numerical orders of convergence of the O-method are usually $\mathcal O(h_{\mathcal M}^2)$
for $\bu$ and $\mathcal O(h_{\mathcal M})$ for the fluxes, provided that $\bu\in H^2$.
If $\bu\in H^{1+\alpha}$ with $\alpha\ge 0$, the orders of convergence
seem to be\cite{AAV06} $\min(2,2\alpha)$ for $\bu$ and $\min(1,\alpha)$ for
its fluxes ($\min(2,\alpha)$ for the fluxes in case of smooth meshes).
It has nonetheless been noticed\cite{AGE10} that, for anisotropy
ratios (the largest eigenvalue of $\Lambda$ divided by the smallest eigenvalue
of $\Lambda$) of order $1000$ or more, the MPFA O-method no longer seems to converge on distorted grids,
due to its loss of coercivity.

L- and G-methods have similar numerical behaviours, but they
seem more stable than the O-method in presence of
strong anisotropy or on irregular meshes used in basin simulation\cite{AAV08,AGE10-2}.

\subsection{Maximum principle for MPFA methods}

When the mesh satisfies the orthogonality condition \eqref{cond-orth2},
MPFA methods are identical to the TPFA scheme and are therefore monotone.
As mentioned, however, such orthogonality conditions are too restrictive
in practice.

For some particular meshes, such as polygonal meshes whose cells are
the union of triangles satisfying the Delaunay condition
(the interaction regions are then triangles),
the O-method is monotone if $\Lambda$
is constant. In the general case, conditions can be found\cite{EIG02} on the
triangle angles and the diffusion tensor to ensure that the
O-method gives rise to an M-matrix, and these conditions can be used
to modify the positions of the mesh vertices in order to try and get
an M-matrix. However, for large anisotropy ratios, such a modification may fail.

In most cases, the L-method displays better monotony properties
than the O-method. The sufficient conditions of Ref. \refcite{NOR07}
(see below) are satisfied by the L-method on a larger class of meshes and tensors
than for the O-method and, even in cases where monotony is violated,
the L-method seems to present much less oscillations than the O-method\cite{AAV08}.

One way to mitigate the problem of large anisotropy in the O-method,
which leads to non-monotony and inaccuracies, is to apply a stretching\cite{AAV98-II}
of the physical space to reduce the anisotropy ratio of $\Lambda$. This
stretching does not seem necessary for regular hexagonal meshes
but mandatory for triangular meshes when the anisotropy ratio is
larger than 10. 

The inaccuracy of the O-method in case of strong anisotropy can also be
reduced by using a variant of the MPFA O-method introduced (in 2D)
separately in Ref. \refcite{CHE08} under the name ``Enriched MPFA O-method'' (EMPFA)
and in Ref. \refcite{EDW08} under the name ``Full pressure support scheme'' (FPS).
This method relaxes the constraints on edge and cell unknowns by adding
vertices unknowns, which gives enough degrees of freedom to
assume the full continuity of the approximation of $\bu$ on the sub-edges (not only at midpoints).
This approximation is taken either piecewise linear (on the triangles $\bxs \V \xcv_K$,
$\xcv_K \V \bx_{\sigma'}$, etc.) or piecewise bilinear (on the subcells
$\V \bx_{\sigma'} \xcv_K\bxs$, $\V \bxs \xcv_L \bx_{\sigma''}$, etc.) and
the new vertex unknown at $\V$ is eliminated by integrating \eqref{base} on a small
domain around $\V$. The monotony (using M-matrix conditions introduced Ref. \refcite{EDW98})
and coercivity of the bilinear variant are analysed for quadrangular meshes
in Ref. \refcite{EDW08} and for triangular meshes in Ref. \refcite{FRI11}.
However, even if the EMPFA/FPS method improves the monotony properties
of the O-method in a number of numerical tests,
it remains unstable (non coercive) in case of strong anisotropy\cite{TRU09}.
According to Ref. \refcite{EDW08,FRI11}, these improved monotony properties
stem from imposing the continuity of the approximation
on whole sub-edges, which prevents the EMPFA/FPS method from displaying
decoupling properties of the O-method shown to be the cause of spurious oscillations.
As mentioned above, the L-method also imposes continuities of full edges and
presents improved monotony characteristics with respect to the O-method
(its extension to 3D meshes moreover appears to be more straightforward than the extension
of cell-centred EMPFA/FPS method). However, to our best knowledge,
numerical or theoretical comparisons of the EMPFA/FPS and L methods
still remain to be done.

A series of interesting results deserves to be mentioned here on the issue
of the monotony of generic 9-point schemes on quadrilateral grids
(which contain the MPFA methods). Sufficient conditions\cite{NOR05,NOR07}
for the monotony of such scheme can be obtained if $\Lambda$
is constant, which provide guidance to generate
meshes on which MPFA methods are monotone, and also show that
7-point methods (such as the L-method) enjoy better monotony
properties in general\cite{AAV08}. These results also prove\cite{KEI09}
that no linear 9-point scheme on generic quadrilateral meshes, which is
exact on linear solutions, can be monotone for any $\Lambda$ (this
has already been noticed, under another form, in Ref. \refcite{KER81}).

\subsection{To summarise: MPFA methods}

The main strengths of MPFA methods are their cell-centred characteristic
and a local computation of the fluxes
(only cell unknowns close to an edge are used in the computation of the flux
across this edge), which lead to acceptable stencils: 9 on 2D quadrilaterals,
27 on 3D hexahedral.
A (small) disadvantage is the necessity to solve local
systems to eliminate the edge/gradient unknowns, which may prove non-invertible
in some cases and therefore require to locally modify the
method\cite{AGE10-2,VOH06}. This however seems to happen relatively
rarely and most numerical tests presented in the literature run without this issue.

A more undesirable characteristic of the MPFA method is their \emph{conditional}
coercivity and monotony. Despite numerous works on the topic, it is not
always obvious to establish \emph{a priori} the range of coercivity or
monotony of an MPFA method on a generic mesh or with a generic diffusion
tensor. As a consequence, unforeseen instabilities and loss of convergence may 
occur.

The question therefore remains to find a FV method which would be
unconditionally coercive and monotone on any type of mesh...

\section{HMM methods}\label{sec:HMM}

Hybrid Mimetic Mixed (HMM) methods are made up of three families of
methods, separately developed in the last ten years or so:
the Hybrid Finite Volume method\cite{EYM10} (HFV),
the Mimetic Finite Difference method\cite{BRE05-I,BRE05-II} (MFD) and the
Mixed Finite Volume method\cite{DRO06} (MFV). It has recently been understood\cite{DRO10}
that all these methods are in fact identical and, therefore, that any analysis
made for one also applies to the other two.

In HMM methods, the main unknowns are cell unknowns $(u_K)_{K\in\mathcal M}$
and edge unknowns $(u_\sigma)_{\sigma\in\mathcal E}$ (approximations of
$(\bu(\bxs))_{\sigma\in\mathcal E}$ where, as in Sec. \ref{sec:MPFA},
$\bxs$ is the centre of gravity of $\sigma$). Of
the three families gathered in HMM methods, MFV methods are the ones
with the most classical FV presentation, involving imposed balance and conservativity
equations \eqref{bf}-\eqref{cf}. Contrary to MPFA methods, edge unknowns are not
eliminated and the computation of the fluxes is made through
local inner products, thus ensuring the coercivity of the scheme.

For given fluxes $F_K=(F_{K,\sigma})_{\sigma\in\mathcal E_K}$ on $\partial K$, we introduce
the vector
\begin{equation}\label{defvmfv}
\gradhmm_K(F_K)=-\frac{1}{|K|}\Lambda_K^{-1}\sum_{\sigma\in\mathcal E_K}F_{K,\sigma}(\bxs-\xcv_K).
\end{equation}
Stokes' formula shows that if $\bu$ is linear in $K$
and $F_{K,\sigma}=-|\sigma|\Lambda_K\nabla\bu_{|K}\cdot \n_{K,\sigma}$, then
$\gradhmm_K(F_K)=\nabla\bu_{|K}$. Hence, $\gradhmm_K(F_K)$ can be considered as a consistent approximation
of $\nabla\bu$ on $K$.
Letting
\begin{equation}
T_K(F_K)=(T_{K,\sigma}(F_K))_{\sigma\in\mathcal E_K}\mbox{ with }
T_{K,\sigma}(F_K)=\frac{1}{|\sigma|}F_{K,\sigma}+\Lambda_K \gradhmm_K(F_K)\cdot\n_{K,\sigma},
\label{defpenmfv}
\end{equation}
the following local inner product is defined
\begin{equation}
[F_K,G_K]_K=|K|\gradhmm_K(F_K)\cdot\Lambda_K\gradhmm_K(G_K)+
T_K(G_K)^T\mathbb{B}_K T_{K}(F_K)
\label{defpslocalmfv}\end{equation}
(where $\mathbb{B}_K$ is a symmetric definite positive matrix)
and the relation between the fluxes and the cell and edge unknowns is
\begin{equation}
\forall G_K=(G_{K,\sigma})_{\sigma\in\mathcal E_K}\in \R^{\mathcal E_K}\,:\;
[F_K,G_K]_K=\sum_{\sigma\in\mathcal E_K}(u_K-u_\sigma)G_{K,\sigma}.
\label{lienpFmfv}\end{equation}

An MFV scheme is defined by \eqref{bf}-\eqref{cf}-\eqref{defpenmfv}-\eqref{defpslocalmfv}-\eqref{lienpFmfv}
for some choices of $(\mathbb{B}_K)_{K\in\mathcal M}$,
with Dirichlet boundary conditions enforced by imposing the value of
$u_\sigma$ if $\sigma\in\mathcal E_{\rm ext}$.
Neumann boundary conditions are as easily considered\cite{CHA07} by imposing 
the value of $F_{K,\sigma}$ for all $\sigma\in\mathcal E_{\rm ext}$.

\begin{remark} For a given edge $\sigma\in\mathcal E_K$, using $G_K(\sigma)=
(\delta_{\sigma,\sigma'})_{\sigma'\in\mathcal E_K}$ ($\delta_{\sigma,\sigma'}$ being
Kronecker's symbol), we can see\cite{DRO10} that
\begin{equation}
u_\sigma-u_K=\gradhmm_K(F_K)\cdot (\bxs-\xcv_K)-T_K(G_K(\sigma))^T\mathbb{B}_KT_K(F_K).
\label{p_increments}\end{equation}
Given that $T_K$ vanishes on exact fluxes of linear functions
and that $\gradhmm_K(F_K)\approx \nabla\bu_{|K}$, \eqref{p_increments} shows that
\eqref{lienpFmfv} is a Taylor expansion
with second order remainder.
\end{remark}

MFD methods are constructed starting from \eqref{lienpFmfv} and looking for
inner products $[\cdot,\cdot]_K$ which satisfy the following consistency condition
(discrete Stokes' formula): for all affine function $q$ and all
$G_K=(G_{K,\sigma})_{\sigma\in\mathcal E_K}\in\R^{\mathcal E_K}$,
\begin{equation}\label{consmfd}
[(\Lambda\nabla q)^I,G]_K+\int_K q(x)(\divmfd G_K)\d x =
\sum_{\sigma\in\mathcal E_K}\frac{1}{|\sigma|}G_{K,\sigma} 
\int_{\sigma}q(x)\d S(x),
\end{equation}
where $((\Lambda\nabla q)^I)_{K,\sigma}=|\sigma|\Lambda_K\nabla q_{|K}\cdot\n_{K,\sigma}$
and $\divmfd G_K=\frac{1}{|K|} \sum_{\sigma\in\mathcal E_K}G_{K,\sigma}$
is the natural discrete divergence of the discrete vector field $G_K$.
{}From the consistency condition \eqref{consmfd}, an algebraic decomposition of
the matrix of $[\cdot,\cdot]_K$(\footnote{i.e. the matrix $\mathbb{M}_K$ 
such that $[F_K,G_K]_K=G_K^T\mathbb{M}_KF_K$.}) can be obtained\cite{BRE05-II} and used
to prove\cite{DRO10} that any inner product satisfying \eqref{consmfd} has the form
\eqref{defpslocalmfv} for some symmetric positive definite $\mathbb{B}_K$.

Relation \eqref{lienpFmfv} can be inverted to express the fluxes in terms
of the cell and edge unknowns and eliminate them. By doing so, we obtain\cite{DRO10}
the HFV scheme. To write down this formulation of the HMM
methods, we introduce for any given vector $u=((u_K)_{K\in\mathcal M},
(u_\sigma)_{\sigma\in\mathcal E})$ the following discrete gradient in $K$:
\begin{equation}
\nabla_K u=\frac{1}{|K|}\sum_{\sigma\in\mathcal E_K}|\sigma|(u_\sigma-u_K)\n_{K,\sigma}.
\label{defgradhfv}\end{equation}
Stokes' formula shows that this gradient is exact if the vector $u$ interpolates
a linear function at $(\xcv_K)_{K\in\mathcal M}$, $(\bxs)_{\sigma\in\mathcal E}$
(it can also be seen\cite{DRO10} that if $u$ and
$F_K$ are related by \eqref{lienpFmfv} then $\nabla_K u=\gradhmm_K(F_K)$).
The function
\begin{equation}
S_K(u)=(S_{K,\sigma}(u))_{\sigma\in\mathcal E_K}\mbox{ with }
S_{K,\sigma}(u)=u_\sigma-u_K-\nabla_K u\cdot(\bxs-\xcv_K)
\label{defSk}\end{equation}
is therefore a first order Taylor expansion, which vanishes
on interpolants of linear functions. The formulation of the HFV method
is then: find $u=((u_K)_{K\in\mathcal M},(u_\sigma)_{\sigma\in\mathcal E})$
(where $u_\sigma$ is fixed by $\udir$ if $\sigma\in\mathcal E_{\rm ext}$)
such that, for any vector $v=((v_K)_{K\in\mathcal M},(v_\sigma)_{\sigma\in\mathcal E})$ with
$v_\sigma=0$ if $\sigma\in\mathcal E_{\rm ext}$,
\begin{equation}
\sum_{K\in\mathcal M}|K|\Lambda_K \nabla_K u\cdot\nabla_K v +\sum_{K\in\mathcal M} S_K(v)^T 
\widetilde{\mathbb{B}}_K S_{K}(u) =\sum_{K\in\mathcal M}v_K\int_K f,
\label{defhfv2}
\end{equation}
where $(\widetilde{\mathbb{B}}_K)_{K\in\mathcal M}$ are symmetric positive
definite matrices (which depend on the matrices $(\mathbb{B}_K)_{K\in\mathcal M}$
in \eqref{lienpFmfv}). This formulation is clearly a discretisation of the
weak formulation \eqref{basew} of \eqref{base}.

\begin{remark} The original MFV, MFD and HFV methods are slightly less general than
the ones presented here. The original MFV method writes \eqref{p_increments} with
a different (stronger) stabilisation, the original MFD method only consider
the case where $\xcv_K$ is the centre of gravity of $K$, and the original HFV
method is only written using diagonal matrices $\widetilde{\mathbb{B}}_K$.
Most of the analysis developed for each of these three methods however extends to the
general HMM method.
\end{remark}

\subsection{Coercivity and convergence of HMM methods}

HMM methods are built on inner products and are therefore unconditionally
coercive (under natural and not very restrictive assumptions on the mesh
regularity). As a consequence and since they are linearly exact,
they enjoy nice stability and convergence properties.
The path of convergence described in Sec. \ref{sec:analysis} has been
successfully applied to HMM methods in Refs. \refcite{DRO06,EYM10}.
Assuming that $\udir=0$ and taking $v=u$ in the discrete variational formulation \eqref{defhfv2}
gives a natural discrete $H^1_0$ norm (the square root of the
left-hand side of the equation), for which one can establish
a Poincar\'e inequality and a discrete Rellich theorem.
The convergence of HMM schemes therefore holds even if $\Lambda$ is discontinuous
and $\bu$ only belongs to $H^1$. For
simplicial meshes, the stabilisation term in \eqref{defpslocalmfv} can
be removed\cite{DRO06} (i.e. $\mathbb{B}_K=0$) without
losing the coercivity, although numerical results are
then slightly less accurate.

Nevertheless, numerical tests\cite{BRE05-II,EYM10} indicate that the choice of
$\mathbb{B}_K$ usually plays little role in the accuracy of the scheme,
provided that this matrix is scaled accordingly to some measure of the eigenvalues of
$\Lambda_K$ (e.g. the trace of this tensor) and that its coercivity properties
incorporate geometric information such as face sizes\cite{DRO10} in case of very distorted
meshes\cite{LIP13-p}. Let us however notice that, in some cases,
$\mathbb{B}_K$ can be selected to ensure the monotony of the
HMM method (see Sec. \ref{sec:HMM-monotone}).

This analysis of HMM method has been extended
to convection-diffusion equations\cite{BEI11}, with various discretisations of the convection
term (centred, upwind, mimetic-based\cite{CAN09}).
General forms of ``automated upwinding'' of the convection, scaled by the local diffusion strength,
are studied in Ref. \refcite{BEI11} and shown to be accurate in all regimes
(diffusion- or convection-dominated).
Numerical experiments also show that much better results are obtained,
in case of strong anisotropy and heterogeneity in a convection-dominated regime, if
the upwinding is made with edge unknowns
rather than cell unknowns (see also Ref. \refcite{DRO10-II} for
the Navier-Stokes equations). This is probably general to many
methods involving edge unknowns, but this would need to be 
theoretically and numerically investigated in a more thorough way.

As HMM methods are based on full gradients reconstructions $\gradhmm_K(F_K)$
or $\nabla_K u$, they are particularly well-suited to non-linear equations
and have been adapted to a number of meaningful models such as
fully non-linear equations of the
Leray-Lions type\cite{DRO06-II} (appearing in particular in models of
non-newtonian fluids), miscible flows in porous media\cite{CHA07}
or the Navier-Stokes equations\cite{DRO09}.
Since the technique in Sec. \ref{sec:analysis} neither
relies on the linearity of the equation nor on the regularity of the solution,
complete convergence analyses of HMM methods for these models are
successfully carried out in these references (along with benchmarking),
under assumptions compatible with applications.

A cell-centred modification (the SUCCES scheme)
of the HMM method, eliminating the edge unknowns by
computing them as convex combinations of cell unknowns,
has been proposed and analysed in Ref. \refcite{EYM10}
for \eqref{base} and in Ref. \refcite{EYM09} for non-linear elliptic equations.
This modification ends up with less unknowns than the HMM method (only cell unknowns) and
is still unconditionally coercive, but it has a larger stencil than
MPFA methods and it displays less accurate numerical results on
grids provoking numerical locking or if $\Lambda$ is discontinuous\cite{EYM08}
(in this latter case, accuracy issues can be mitigated by retaining edge unknowns at
the discontinuities, giving rise to the SUSHI scheme).

When $(\xcv_K)_{K\in\mathcal M}$ are the centres of gravity of the cells, 
HMM methods are the original (edge-based) MFD methods and
all results on these methods apply to HMM methods, for example:
convergence rates for smooth data
and super-convergence of $u$ if a proper lifting of the numerical
fluxes exists\cite{BRE05-I,BRE07}, \emph{a posteriori} estimators
usable for mesh refinement\cite{BEI08,BEI08-II},
higher order methods designed to recover optimal orders of convergence
on the fluxes\cite{GYR08,BEI08-I,BEI09}, or extension to non-planar
faces\cite{BRE06,BRE07,LIP06}.
We will not delve into more details here and we refer to Ref. \refcite{LIP13}
for a comprehensive review of MFD methods. One
open issue however seems interesting to mention regarding the extensions of MFD methods
which introduce additional flux unknowns (higher order methods or methods for
non-planar faces). These methods are based on the construction of
local scalar products satisfying a generalisation of the consistency relation \eqref{consmfd}
on the expanded flux space.
Algebraic decomposition of these scalar product matrices are known\cite{BEI08-I,BRE07},
but the question remains open to find expression of these products
purely based on geometrical quantities such as in \eqref{defpslocalmfv}.
This would in particular eliminate the need to solve local algebraic problems to construct them.

\begin{remark}[Mixing MPFA and HMM ideas]\label{rem:mix}
In Refs. \refcite{AGE09,EYM12}, the sub-cells flux continuity of the MPFA methods
is combined with the gradient and stabilisation \eqref{defgradhfv}-\eqref{defSk}
of HMM methods (on the same sub-cells, by introducing half-edge unknowns) to construct an unconditionally coercive and convergent
scheme. If the mesh and diffusion tensors are not too skewed,
the sub-cells can be defined using particular harmonic edge points (instead of
$\bxs$), where the solution can be interpolated using only the two
neighbouring cell values. In this case, the half-edge unknowns can
be eliminated vertex by vertex, as in the O-method, and a 9-point stencil cell-centred
scheme is recovered on quadrilateral meshes.

Another mixing of MPFA and HMM ideas can be found in the method presented in Ref. \refcite{LIP09}.
This method uses, as the MPFA O-method,
additional face unknowns (as many on $\sigma$ as the number of vertices of $\sigma$)
but constructs local ``scalar products'' in each sub-cell around
a given vertex, trying to satisfy the local consistency conditions
\eqref{consmfd}. Except on simplicial meshes, construction of
such consistent coercive scalar products is not theoretically proved,
but when they exists their block structure around each vertex allows one,
as in the O-method, to eliminate the face unknowns
and obtain a coercive method with the same stencil as the O-method.

\end{remark}

\begin{remark}[Mixing HMM, MPFA and dG ideas]
Ref. \refcite{DIP12} proposes a scheme which mixes
HMM, MPFA and dG ideas. This method consists in
constructing a finite-dimensional subspace $V_h$ of piecewise
affine functions, whose gradient in each cell is
given by \eqref{defgradhfv} in which the edge unknowns are
computed from cell unknowns using the elimination technique of the MPFA L-method.
This space $V_h$ is then used in a Finite-Element like discretisation of
\eqref{basew} with a bilinear form including jumps penalisations
as in dG methods.

If the edges unknowns are not eliminated then numerical fluxes can
be found\cite{DIP13} such that this scheme satisfies the balance and conservativity
equations \eqref{bf}-\eqref{cf}.
\end{remark}

\subsection{Maximum principle for HMM methods}\label{sec:HMM-monotone}

HMM methods are usually not monotone, even on parallelogram meshes and for
constant $\Lambda$. In simple cases, one can obtain
necessary and/or sufficient conditions on the diffusion tensor and the mesh for the existence
(i.e. a choice of $\mathbb{B}_K$) of a monotone HMM method\cite{LIP11,LIP11-II}.
The idea is to hybridise the method (i.e. eliminate the cell unknowns, see
Sec. \ref{sec:HMM-summary}) and to analyse if the corresponding matrix
is an M-matrix and if the corresponding right-hand side is non-negative
whenever $f\ge 0$.

For simplicial meshes, a necessary and sufficient condition of monotony of
any HMM method is that $\Lambda_K\mathbf{n}_{K,\sigma}\cdot\mathbf{n}_{K,\sigma'}<0$
for all $K\in\mathcal M$ and all $\sigma\not=\sigma'\in\mathcal E_K$
(if $\Lambda$ is isotropic, this comes down to imposing that all angles
of the simplicial meshes are less that $\pi/2$). \emph{Necessary} monotony conditions
can be written for meshes made of parallelograms or parallelepipeds, which turn
out to be identical to the conditions in 2D for
9-point cell-centred schemes\cite{NOR07}. These conditions give insights on how to
construct, using the algebraic point of view of MFD methods,
the matrices of the local scalar products $[\cdot,\cdot]_K$ in \eqref{defpslocalmfv},
but remain to be translated into \emph{geometric} constructions
of proper $\mathbb{B}_K$ matrices. Although similar conditions can also
be written for other types of meshes, such as locally refined rectangular
meshes\cite{LIP11}, a more thorough analysis remains to be done to find necessary
and/or sufficient monotony conditions for HMM methods on generic meshes. 
Ref. \refcite{LIP11-II} suggests, in the absence of such an analysis, to use a heuristic
based on constructing $\mathbb{B}_K$ by solving local optimisation problems
which penalise the scalar products $[\cdot,\cdot]_K$ whose matrix is not
an M-matrix.

\subsection{Coercivity vs. Monotony vs. Accuracy}

\begin{figure}[h!]
\begin{center}
\begin{tabular}{ccc}
\includegraphics[width=0.28\linewidth]{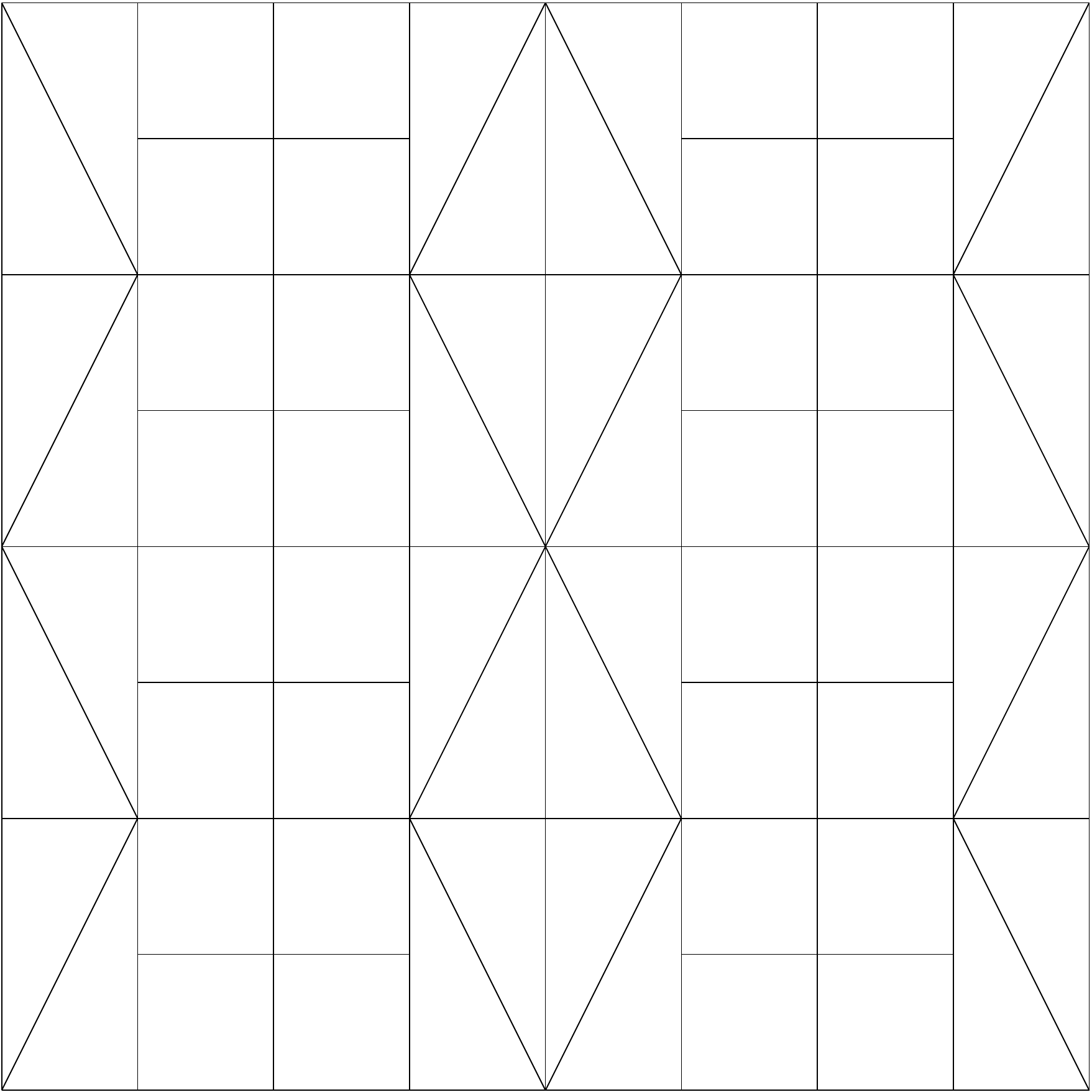}&
\includegraphics[width=0.28\linewidth]{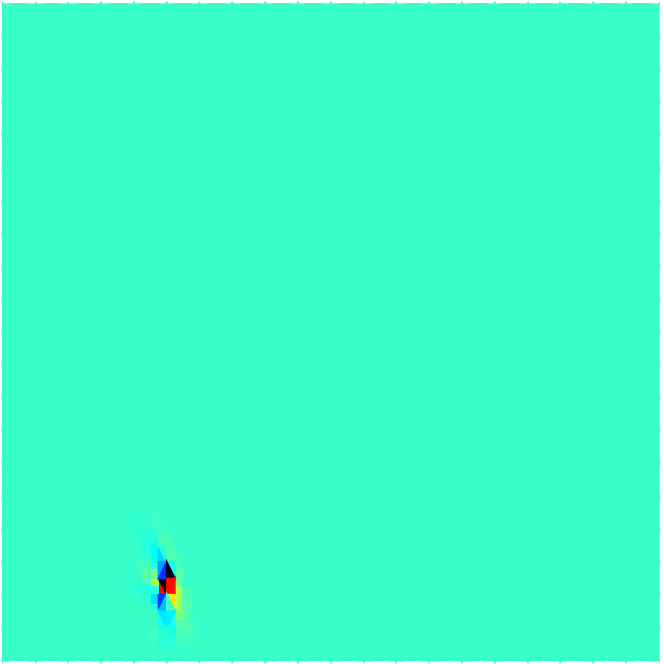}&
\includegraphics[width=0.28\linewidth]{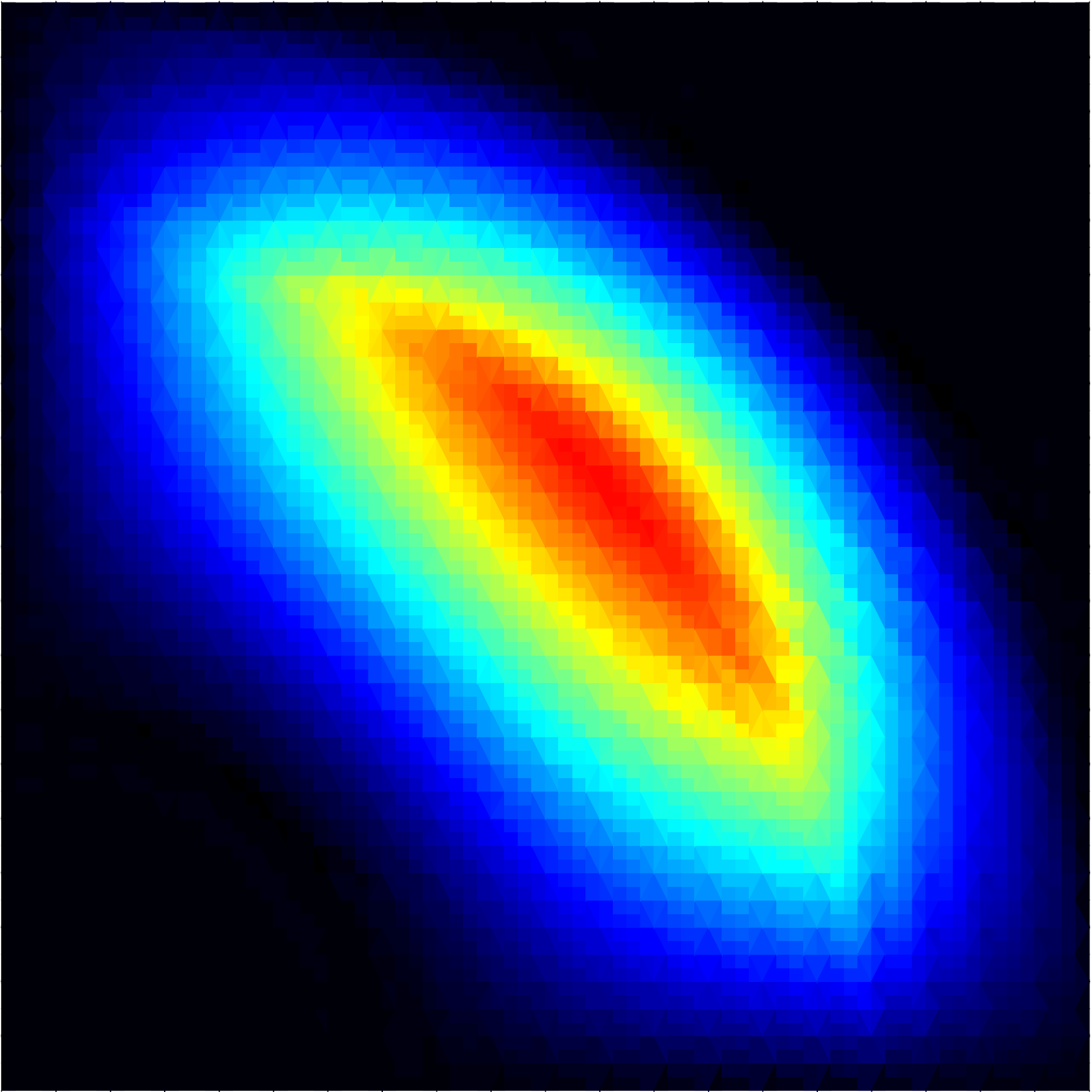}\\
Mesh pattern&G-scheme solution,&HMM solution,\\
(mesh=$10\times 10$ reproduction&$\min u=-9\times 10^{240}$&$\min u=-7.9\times 10^{-3}$\\
of this pattern)&
$\max u=7\times 10^{240}$&$\max u=0.52$
\end{tabular}
\caption{\label{fig:Gexplose}Explosion of a non-coercive method applied to a
transient problem.}
\end{center}
\end{figure}

If a scheme's matrix has negative eigenvalues,
any negative mode will be amplified when the scheme is applied
to a transient equation, thus provoking the explosion of the solution.
Fig. \ref{fig:Gexplose} illustrates this phenomenon when a (non-coercive) G-scheme
and a time-implicit discretisation (involving 150 time steps)
is applied with $\Omega=(0,1)^2$ and final time $T=0.1$ to
$\partial_t \bu -\div(\Lambda\nabla\bu)=0$ with $\udir=0$ and
\[
\Lambda(x,y)=\frac{1}{x^2+y^2}\left(\begin{array}{c@{\quad}c} 10^{-3}x^2+y^2 & (10^{-3}-1)xy\\
							(10^{-3}-1)xy & x^2+10^{-3}y^2\end{array}\right)\,,\;\;
u(0,\cdot)=\left\{\begin{array}{ll} 1&\mbox{ on $(\frac{1}{4},\frac{3}{4})^2$},\\
	0&\mbox{ elsewhere.}\end{array}\right.
\]
The coercivity of a scheme does not only ensure that it converges as the mesh is refined,
but also that it does not explode in transient cases as shown for the HMM method in Fig. \ref{fig:Gexplose}
(the HMM solution is quite close to the expected solution in this test case).

\begin{figure}[h!]
\begin{center}
\begin{tabular}{ccc}
\includegraphics[width=0.28\linewidth]{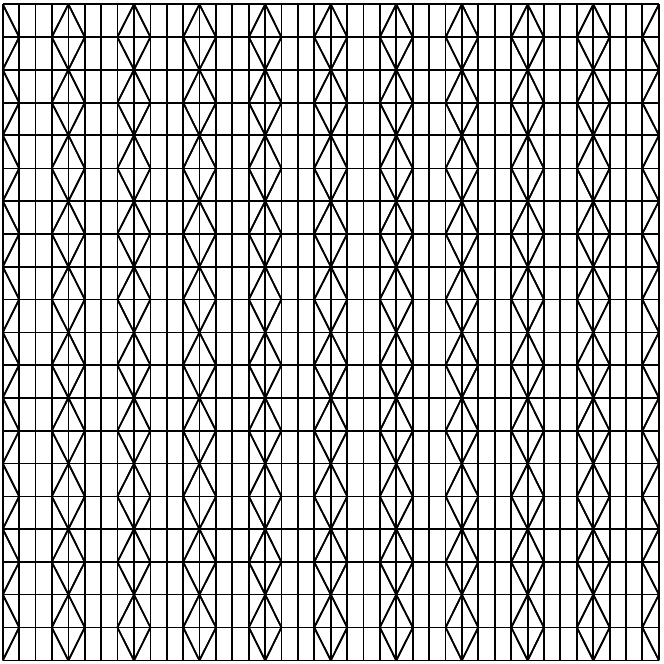}&
\includegraphics[width=0.28\linewidth]{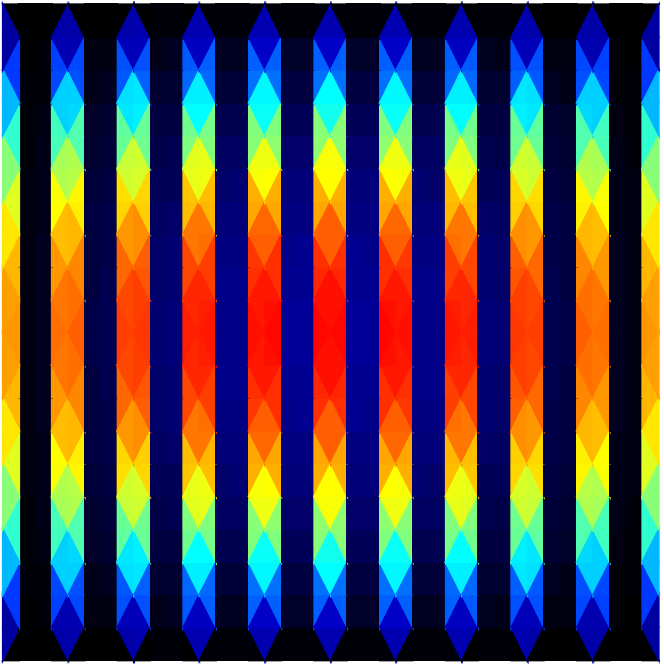}&
\includegraphics[width=0.28\linewidth]{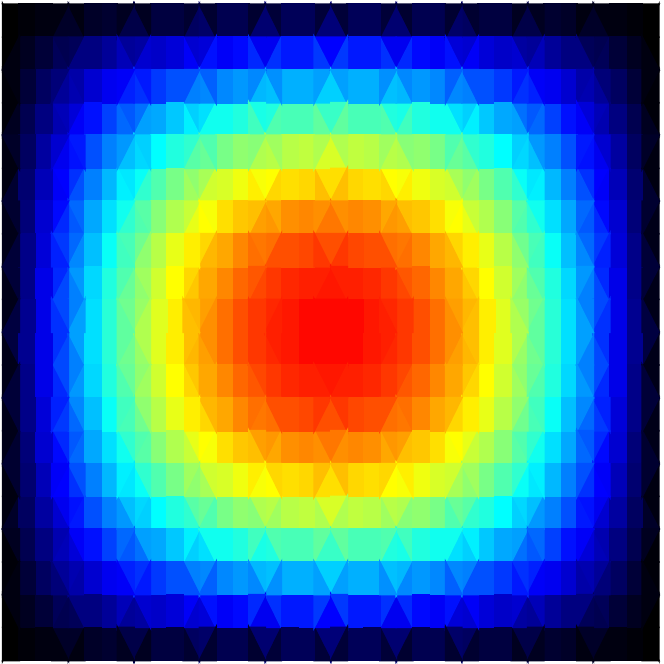}\\
Mesh&HMM&G-scheme
\end{tabular}
\caption{\label{fig:hmmG}Numerical test with strong anisotropy ratio.
The G-scheme is not coercive
in this test case.}
\end{center}
\end{figure}

The convergence insured by the coercivity of a method however does not
mean that it is always accurate (only that it is accurate
as the mesh size tends to $0$). For instance, the unconditionally coercive
HMM method may display very bad numerical behaviour in presence
of strong misalignment between the grid directions and the
principal directions of diffusion.
In Fig. \ref{fig:hmmG}, we present the numerical
solutions produced by an HMM method and the G-scheme
for the constant diagonal tensor $\Lambda={\rm diag}(10^4,1)$
and the exact solution $\bu(x,y)=x(1-x)y(1-y)$.
The strong oscillations displayed by the HMM
method in this example are probably due to its lack of monotony
properties and to its non-local computations of
the numerical fluxes ($F_{K,\sigma}$ is expressed in term of \emph{all}
the edge unknowns around $K$, not just unknowns around $\sigma$).
Although it can be checked that the G-scheme is \emph{not} coercive
(and therefore not monotone) on this test case, its
local computation of the fluxes prevents its solution from presenting
spurious oscillations, and therefore seems to improve its 
``apparent'' monotony properties.

\subsection{To summarise: HMM methods}\label{sec:HMM-summary}

The strength of HMM methods is their unconditional coercivity,
on any mesh and for any diffusion tensor. This is achieved
at the cost of a larger number of unknowns (cell and edge
unknowns) than in MPFA methods, but hybridisation techniques can be applied as in Mixed Finite Element
methods to locally eliminate the cell unknowns and retain only the edge
unknowns. This unconditional coercivity ensures the robustness of
HMM methods (no explosion for transient equations) and provides the means
for full convergence analyses for a vast range of different complex models, involving
non-linearities and non-smooth data and solutions.

HMM methods are however not always monotone and, despite the
large freedom in their construction (through the choice
of the matrices $\mathbb{B}_K$), the analysis of their monotony range
is to date very limited. Another weakness is their relative non-local
computation of the fluxes, as $F_{K,\sigma}$ depends on all
edge unknowns around $K$. Because of this, they may present inaccurate results
on coarse meshes in presence of strong anistropy -- although their
unconditional coercivity ensures that, as the mesh is refined, the approximate
solution converges to the exact solution.

The question still remains to find a FV method which would be
unconditionally coercive and monotone on any type of mesh...

\section{DDFV methods}\label{sec:DDFV}

Discrete Duality Finite Volume (DDFV) methods have been introduced 
around the early 2000's\cite{HER98,HER00,HER03}, but have been mostly
studied after 2005\cite{DOM05,AND07,BOY08-II}.
The basic idea of DDFV methods in 2D is a bit similar to
MPFA methods and also draws some inspiration from Ref. \refcite{COU99}.
The initial remark is that the two values $u_K$ and $u_L$ around $\sigma$
only give an approximation of the local gradient in the direction
$(\xcv_K\xcv_L)$ and are therefore insufficient to obtain an
expression of the whole gradient around $\sigma$ (when the orthogonality condition
\eqref{cond-orth2} does not hold,
the whole gradient is required to compute an approximate flux
${F}_{K,\sigma}$).
So, as in MPFA methods, DDFV methods introduce new unknowns to get
an approximation of the gradient in another direction than
$(\xcv_K\xcv_L)$. Using these approximate projections
of the gradient on two independent directions, an approximation of the whole gradient
can be reconstructed in a similar way as \eqref{lin-grad}
defines the gradient \eqref{eq-grad} in MPFA methods.

\begin{figure}[h!]
\begin{center}
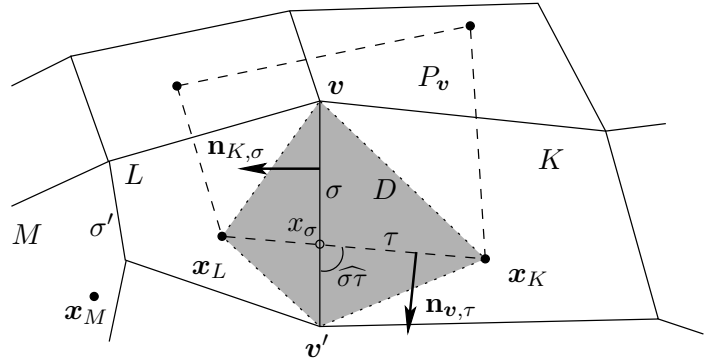
\caption{\label{fig:ddfv}DDFV primal meshes (continuous lines: $K$, $L$, $M$),
dual meshes (dashed lines: $P_\V$) and diamonds (filled: $D$).
$\n_{K,\sigma}$ and $\n_{\V,\tau}=$ unit normals
to $\sigma=[\V,\V']$ and $\tau=[\xcv_K,\xcv_L]$.}
\end{center}
\end{figure}

The additional unknowns of DDFV methods are located at the vertices of the mesh
(we denote by $\mathcal V$ the set of vertices and we refer to
Fig. \ref{fig:ddfv} for notations).
{}From the cell $(u_K)_{K\in\mathcal M}$ and
vertex $(u_\V)_{\V\in\mathcal V}$ unknowns and since $\vect{\V\V'}$ and
$\vect{\xcv_K\xcv_L}$ are linearly independent,
a constant approximate gradient $\nabla_D u$ can be computed
on the diamond $D:={\rm co}(\sigma\cup\{\xcv_K\})
\cup{\rm co}(\sigma\cup\{\xcv_L\})$(\footnote{``${\rm co}$'' denotes the convex hull.
Note that the diamond $D$ may be non-convex (this is 
the case for the diamond around $\sigma'$ in Fig. \ref{fig:ddfv}).})
by imposing $\nabla_D u\cdot (\xcv_K-\xcv_L)=u_K-u_L$ and 
$\nabla_D u\cdot (\V-\V')=u_\V-u_{\V'}$, which leads to\cite{DOM05,AND07}
\begin{equation}\label{ddfv-eqgrad}
\begin{array}{lcl}
\dsp\nabla_D u &=&\dsp \frac{1}{\sin(\widehat{\sigma\tau})}
\left(\frac{u_L-u_K}{\d(\xcv_K,\xcv_L)}\n_{K,\sigma}+\frac{u_{\V'}-u_\V}{\d(\V,\V')}\n_{\V,\tau}\right)\\
&=&\dsp \frac{1}{2|D|}\left((u_L-u_K)\d(\V,\V')\n_{K,\sigma}+(u_{\V'}-u_\V)\d(\xcv_K,\xcv_L)\n_{\V,\tau}\right),
\end{array}
\end{equation}
where $\widehat{\sigma\tau}$ is the angle between the straight
lines $(\xcv_K\xcv_L)$ and $(\V\V')$ and $|D|$ is the area of $D$.
One can then compute an approximate flux
through $\sigma$:
\begin{equation}\label{ddfv-flux1}
F_{K,\sigma}=-|\sigma|\Lambda_D \nabla_D u\cdot\n_{K,\sigma},
\end{equation}
where $\Lambda_D$ is the mean value of $\Lambda$ on $D$. The balance
equations on each cell \eqref{bf} then give as many equations as 
the number of cell unknowns. To close the system, it remains to find as many
equations as the number of vertex unknowns, which is simply done by
writing the balance equation on new cells (``dual cells'') constructed around vertices.
A natural choice\cite{DOM05,AND07,BOY08-II,HER00} for the dual
cell around $\V$ is the polygon $P_\V$ which has all
the cell points $\xcv_K,\xcv_L,\ldots$ around $\V$ as vertices (in dotted
lines in Fig. \ref{fig:ddfv}). The flux through the edge $\tau=[\xcv_K,\xcv_L]$
of $P_\V$ can be computed using the gradient on $D$:
\begin{equation}\label{ddfv-flux2}
F_{\V,\tau}=-|\tau|\Lambda_D \nabla_D u\cdot\n_{\V,\tau}
\end{equation}
and the balance of these fluxes around a vertex $\V$ reads
\begin{equation}\label{ddfv-bf2}
\sum_{\tau\in \mathcal E_{P_\V}}F_{\V,\tau}=\int_{P_\V} f(x)\d x.
\end{equation}
where $\mathcal E_{P_\V}$ is the set of all edges of $P_\V$.
These balance equations around each vertex complete the set
of equations which define the DDFV method, that is
\eqref{bf}-\eqref{ddfv-eqgrad}-\eqref{ddfv-flux1}-\eqref{ddfv-flux2}-\eqref{ddfv-bf2}.
Note that the flux conservativity across primal $\sigma$ and dual $\tau$
edges are naturally satisfied by \eqref{ddfv-flux1} and \eqref{ddfv-flux2}.

\begin{remark}
Dirichlet or Neumann boundary conditions are handled seamlessly. The diamond
around a boundary edge $\sigma\in\mathcal E_K\cap\mathcal E_{\rm ext}$
is only made of the triangle ${\rm co}(\sigma\cup\{\xcv_K\})$, and the gradient
on $D$ is constructed by replacing $\xcv_L$ with a point $\xcv_\sigma\in\sigma$
(which is also used to define the dual cell around $\V$) and $u_L$
with some unknown $u_\sigma$. 
Dirichlet boundary conditions then fix $(u_\sigma)_{\sigma\in\mathcal E_{\rm ext}}$ and
$(u_\V)_{\V\in \mathcal V\cap\partial\Omega}$ using the values of $\udir$, and
\eqref{ddfv-bf2} is not written for boundary vertices\cite{AND07,DOM05}.
Neumann boundary conditions simply impose the value of $F_{K,\sigma}$,
and \eqref{ddfv-bf2} is written for all vertices\cite{DOM05}.
\end{remark}

The preceding construction is valid if all dual cells $(P_\V)_{\V\in\mathcal V}$
have disjoint interiors and, therefore, form a partition $\Omega$. It may happen for
peculiar meshes that the preceding construction of $P_\V$ leads to overlapping
dual cells. In this case, the scheme must be modified and a possible choice\cite{HER03} is
to take for $P_\V$ the interaction region around $\V$ from the MPFA
methods (see Fig. \ref{fig:mpfa}).

If $\Lambda$ is discontinuous across $\sigma$, the usage of
its mean value on $D$ in \eqref{ddfv-flux1} and \eqref{ddfv-flux2} may lead to a loss of accuracy.
In case this case, and still assuming that $\Lambda$ is
constant on each (primal) cell $K\in\mathcal M$, the DDFV scheme can
be modified\cite{HER03} by introducing an unknown $u_\sigma$
at the point $\{\xcv_\sigma\}=\sigma\cap (\xcv_K\xcv_L)$ (or $\bxs$ if $D$ is not convex
and $P_\V$ is the same interaction region as in MPFA methods),
using it to compute constant gradients in each half-diamond
$D\cap K$ and $D\cap L$ and then eliminating it thanks to the flux conservativity \eqref{cf}
through primal edges. Since there is no jump of $\Lambda$ 
through $\tau=[\xcv_K,\xcv_L]$, the conservativity through this
dual edge is ensured as the sub-fluxes through $[\xcv_K,\xcv_\sigma]$ and $[\xcv_\sigma,\xcv_L]$
use the same values of $\Lambda$ on each side of $\tau$ (respectively $\Lambda_K$
and $\Lambda_L$) and the same gradient on each half diamond.
If $\Lambda$ is also discontinuous across dual edges (which is not standard
in reservoir engineering), a further
modification of the DDFV method has been proposed in Ref. \refcite{BOY08-II}.
This ``m-DDFV'' method uses local gradients which are
constant in quarters of diamonds. Four new
unknowns need to be introduced in each diamond, and are then eliminated by imposing (as in MPFA
methods) flux conservativity equations through the diamond diagonals.

Although this presentation of DDFV methods clearly shows that they are
based on FV principle (flux conservativity and balance), it does not
explain the name ``Discrete Duality Finite Volume''. DDFV methods
can be re-cast using discrete gradient and divergence operators, in such
a way that the Green-Stokes duality formula holds at the discrete
level\cite{DOM05,AND07,BOY08-II}. The gradient operator, already
defined, takes cell and vertex values (assumed to represent
piecewise constant functions in primal and dual cells) and constructs a piecewise
constant gradient on the diamonds. The divergence operator
takes a piecewise constant vector field $(\xi_D)_D$ on diamonds and defines
its divergence as piecewise constant functions on primal and dual cells
by writing the flux balances \eqref{bf} and \eqref{ddfv-bf2}
with $F_{K,\sigma}=|\sigma|\xi_D\cdot\n_{K,\sigma}$ and $F_{\V,\tau}=|\tau|\xi_D 
\cdot\n_{\V,\tau}$. Under this form, DDFV methods are based on similar
principles as MFD methods, which aim at satisfying the discrete
Green-Stokes formula \eqref{consmfd}. They are
different methods 
but DDFV methods can be re-cast in a framework similar to MFD methods\cite{COU10}.

Generalisation of DDFV methods to 3D is
based on similar ideas as in the 2D case, but requires quite heavy notations to be
properly defined. Two essentially different 3D generalisations exist:
methods using Cell and Vertex unknowns (hence dubbed CeVe-DDFV)
and methods relying on Cell, Vertex, Faces and Edges unknowns
(called CeVeFE-DDFV).
Refs. \refcite{HER09,COU09,AND10} design CeVe-DDFV methods by
reconstructing a piecewise constant gradient from
its projection on $\vect{\xcv_K\xcv_L}$ computed using $u_K$ and $u_L$,
and its projection on the plane generated by $\sigma$
computed using the values on the vertices of $\sigma$.
Linearly exact formulas can be found for this projected gradient\cite{AND12}
but the discrete Poincar\'e inequality (crucial to Step \textsf{(C1)}
in Sec. \ref{sec:analysis}) only seems provable
when all faces $\sigma$ are triangles(\footnote{Or on cartesian grids\cite{AND13}.})
and the CeVe-DDFV method is therefore not coercive on generic meshes.
Refs. \refcite{COU11,COU11-II} propose a CeVeFE-DDFV method
with a local gradient computed from its projection on $\vect{\xcv_K\xcv_L}$ and
$\vect{\V\V'}$ (as in 2D) and on a third face-edge direction.
A third mesh is built around each face and edge centres to
obtain additional balance equations for the new face and edge unknowns.
This CeVeFE-DDFV method is coercive on any mesh, but at the cost
of additional unknowns with respect to the CeVe-DDFV method.

\subsection{Coercivity and convergence of DDFV methods}

Because DDFV methods are based on discrete gradient and divergence
operators which reproduce, as MFD methods, the Green-Stokes formula,
discrete $H^1_0$ estimates can be obtained by mimicking the continuous
integration by parts \eqref{ipp}, provided that the discrete
Poincar\'e inequality holds. This is the case in 2D,
for the CeVeFE-DDFV 3D method or for the CeVe-DDFV 3D method on
meshes with triangular faces. In these cases,
DDFV methods are coercive and, being linearly
exact, they enjoy the corresponding stability and convergence properties.

The technique outlined in Sec. \ref{sec:analysis} has been applied\cite{AND07}
to prove the convergence, without additional
regularity assumption on the data or the solution, of
the 2D DDFV method using the mean values $\Lambda_D$
as in \eqref{ddfv-flux1}-\eqref{ddfv-flux2} (Ref. \refcite{AND07}
provides in fact a convergence analysis for a non-linear equation,
which contains \eqref{base} as a particular case). An
$\mathcal O(h_{\mathcal M})$ error estimate for $u$ and the discrete gradient are
also established if $\Lambda$ is Lipschitz-continuous and
 $\bu\in H^2$ (this estimate was known\cite{DOM05} for $\Lambda={\rm Id}$).

Concerning the m-DDFV method\cite{HER03,BOY08-II}, 
an $\mathcal O(h_{\mathcal M})$ error estimate for $\bu$ and its
gradient has been proved in Ref. \refcite{BOY08-II} (also for a non-linear
version of \eqref{base}), provided that $\bu$ is $H^2$ on each half- or quarter-diamond.
This regularity assumption does not seem always satisfied
(in particular if $\Omega$ or some cells around discontinuities
of $\Lambda$ are not convex), but the path described in
Sec. \ref{sec:analysis} could also be applied to the m-DDFV method.

An $\mathcal O(h_{\mathcal M})$ error estimate on $\bu$ has
been obtained in Ref. \refcite{COU11} for the 3D CeVeFE-DDFV method,
under the assumptions that $\Lambda$ is Lipschitz-continuous and
that $\bu\in H^2(\Omega)$. We can however notice 
that this CeVeFE-DDFV method (as well as the 2D DDFV scheme)
is a Gradient Scheme\cite{EYM12} and, therefore, that
its convergence without regularity assumptions, for \eqref{base} as well as
non-linear and non-local equations,
follows from the general convergence analysis of Gradient Schemes\cite{EYM12,DRO12}.

As HMM methods, DDFV methods have been adapted to more complex models
than \eqref{base}: non-linear elliptic equations\cite{AND07,BOY08-II,COU11},
stationary and transient convection-diffusion equations\cite{COU10,HER12},
the cardiac bidomain model\cite{AND11}, div-curl problems\cite{DEL07,HERM08},
degenerate hyperbolic-parabolic problems\cite{AND10} (with assumptions on the mesh,
see Sec. \ref{sec:ddfv-monotone}), the linear Stokes equations with
varying viscosity\cite{KRE11,KRE12}, semiconductor models\cite{CHA09}
and the Peaceman model\cite{CHA13}.
The convergence analysis of DDFV methods is carried out (sometimes under regularity
assumptions) for all these models except the last two.
Analysis tools for the 3D CeVe-DDFV method are presented in
Ref. \refcite{AND13,AND12} and used to study its convergence
for transient non-linear equations or systems.

\subsection{Maximum principle for DDFV methods}\label{sec:ddfv-monotone}

On meshes satisfying the orthogonality conditions \eqref{cond-orth} or
\eqref{cond-orth2}, DDFV methods for \eqref{base} are identical to
two TPFA schemes\cite{DOM05} (one on each primal and dual mesh),
and are therefore monotone. This monotony
under orthogonality conditions on the mesh is used in
Ref. \refcite{AND10} to study DDFV discretisations of degenerate hyperbolic-parabolic
equations, and in particular to establish discrete entropy inequalities
on approximate solutions. Study of the
monotony of DDFV methods on generic meshes however remains to be done.

\subsection{To summarise: DDFV methods}

As HMM methods, the main strength of DDFV methods is their unconditional
coercivity (with some caveats for 3D methods, see above),
which ensures their robustness and allows one to adapt them and
analyse their convergence for a number of models. Another very practical
property for the analysis of DDFV methods is their
discrete duality property (existence of discrete gradient and divergence operators satisfying
Stokes' formula), which is also shared by HMM methods.
An advantage of DDFV methods over HMM methods is perhaps their
more local computation of the fluxes
($F_{K,\sigma}$ is expressed in terms of unknowns localised around the 
edge $\sigma$, whereas in HMM methods this flux requires all
edge unknowns around $K$).

A relative weakness of DDFV methods is their intricacy in 3D.
The heavy and numerous notations required for the definitions of 3D
DDFV methods probably makes them difficult to adopt by non-specialists
and complexifies their analysis. In particular, establishing the discrete
duality formula is far from obvious. Once passed these complicated
notations, however, implementation of 3D DDFV methods is not particularly difficult.
The lack of monotony studies for DDFV methods is also
a gap in the literature, which would probably need to be filled
to get a better understanding on the possible applicability of these methods
to multi-phase flow models.

And so our quest for an unconditionally coercive and monotone FV
method on any mesh continues...

\section{Monotone and Minimum-Maximum preserving (MMP) methods}\label{sec:LMP}

Previously cited results\cite{NOR05,NOR07,KEI09,KER81} show that
no \emph{linear} 9-point scheme on quadrangular meshes, exact for linear functions (i.e.
of formal order 2), can be monotone on any distorted mesh or for any
diffusion tensor. Some constraints must be relaxed...
One choice is to allow for larger stencils (see Ref. \refcite{LEP09-II} for a Finite Difference scheme).
For Finite Volume methods, the most common choice appears to be a relaxation
of the \emph{linearity} of the scheme and the construction of \emph{non-linear}
``monotone'' approximations of the linear equation \eqref{base}.
The obvious trade-of is that computing the solution to the scheme
will be more complex, requiring Picard or Newton iterations,
which may create computational issues (such as the choice of stopping
criteria). Also, the monotony, conservativity and/or consistency may
only be achieved for the genuine solution to the non-linear scheme,
not at each iteration of these algorithms\cite{LEP09}.

Contrary to MPFA, HMM or DDFV methods, schemes presenting
discrete minimum-maximum principles do not form a well defined family of methods but
are rather schemes constructed using similar ideas and trying to achieve
the discrete minimum principle \eqref{disc-minpple} or
the discrete minimum-maximum principle \eqref{disc-minmaxpple}.
As we are considering non-linear schemes, these two principles
are not equivalent and we should make sure that we clearly separate
both. Schemes satisfying \eqref{disc-minpple} will be called \emph{monotone},
as a commonly used but somewhat misguided extension of the vocabulary
used for linear schemes(\footnote{Indeed, ``monotone'' non-linear methods
do not necessarily preserve orders of boundary conditions or of initial condition for time-dependent problems.
They merely provide solutions which remain non-negative when the
boundary/initial conditions are non-negative.}), whereas schemes which
satisfy \eqref{disc-minmaxpple} will be called \emph{minimum-maximum
preserving} (MMP) schemes.

A widespread idea to obtain a monotone or MMP scheme
is to compute two \emph{linear} fluxes $F^1_{K,\sigma}$
and $F^2_{K,\sigma}$ for each interior edge and to define $F_{K,\sigma}$
as a convex combination of
$F^1_{K,\sigma}$ and $F^2_{K,\sigma}$ with coefficients depending upon the unknown $u$:
\begin{equation}\label{combconv}
\begin{array}{l}
\dsp F_{K,\sigma}=\mu^1_{K,\sigma}(u)F^1_{K,\sigma} + \mu^2_{K,\sigma}(u)F^2_{K,\sigma}\\[0.5em]
\dsp \mbox{with
$\mu^1_{K,\sigma}(u)\ge 0$, $\mu^2_{K,\sigma}(u)\ge 0$ and $\mu^1_{K,\sigma}(u)+
\mu^2_{K,\sigma}(u)=1$}.
\end{array}
\end{equation}
The methods we consider here are cell-centred, but
the definition of $F^1_{K,\sigma}$ and $F^2_{K,\sigma}$ may require
to introduce additional unknowns (e.g. vertex, edge or other unknowns).
These unknowns are then eliminated, classically by expressing them as convex combinations
of cell unknowns.
The coefficients $\mu^1_{K,\sigma}(u)$ and $\mu^2_{K,\sigma}(u)$ are chosen to
eliminate the ``bad'' parts of $F^1_{K,\sigma}$ and $F^2_{K,\sigma}$,
responsible for the possible loss of monotony.

\subsection{Non-linear ``2pt-fluxes'': monotone schemes}\label{sec:nl2pt}

The TPFA scheme is monotone thanks to its ``2pt-flux'' structure. This suggests to
try and build monotone methods on generic meshes by computing
$F_{K,\sigma}$ with a ``2pt'' formula, involving apparently only $u_K$ and $u_L$
but with coefficients depending on all
cell unknowns and boundary values $U=((u_M)_{M\in\mathcal M},(u_\sigma)_{\sigma\in\mathcal E_{\rm ext}})$
(same notation as in Sec. \ref{sec:monotony}). Indeed, assume that $F_{K,\sigma}$
is written
\begin{equation}\label{nonlin-2pt}
F_{K,\sigma}=\alpha_{K,L}(U) u_K - \beta_{K,L}(U) u_L
\quad\mbox{with $\alpha_{K,L}(U)> 0$ and $\beta_{K,L}(U)> 0$}
\end{equation}
(where $L$ is the cell on the other side of $\sigma\in\mathcal E_K\cap
\mathcal E_{\rm int}$ and $L=\sigma$ whenever $\sigma\in\mathcal E_K\cap \mathcal E_{\rm ext}$).
Then the conservativity relation \eqref{cf} imposes,
assuming that it must be satisfied for any value of $U$,
\begin{equation}\label{nonlin-cf}
\alpha_{K,L}(U)=\beta_{L,K}(U)\mbox{ for any neighbour cells $K$ and $L$}.
\end{equation}
The scheme \eqref{bf} can then be recast as 
\begin{equation}\label{nl-recast}
A(U)(u_K)_{K\in\mathcal M}=(B_K(U))_{K\in\mathcal M},
\end{equation}
where $B_K(U)=\int_K f(x)\d x + \sum_{\sigma\in\mathcal E_{\rm ext}\cap
\mathcal E_K}\beta_{K,\sigma}(U)u_\sigma$
and the matrix $A(U)$  has
(i) diagonal coefficients $A_{K,K}(U)=\sum_{M}\alpha_{K,M}(U)>0$
(the sum being on all $M$ neighbour cells or edges of $K$),
(ii) extra-diagonal coefficients $A_{K,L}(U)=-\beta_{K,L}(U)<0$ if
$K$, $L$ are neighbour cells, $A_{K,L}(U)=0$ otherwise,
and (iii) is diagonally dominant by column (strictly for columns
$L$ such that $\mathcal E_{\rm ext}\cap\mathcal E_L\not=\emptyset$)
thanks to \eqref{nonlin-cf}. The graph of $A(U)$ is also connected and
(cf. Sec. \ref{sec:monotony}) $A(U)^{-1}$ therefore
has non-negative coefficients, which means that
the scheme \eqref{nl-recast} satisfies \eqref{disc-minpple}.

\subsubsection{Triangular meshes}

\begin{figure}[h!]
\begin{center}
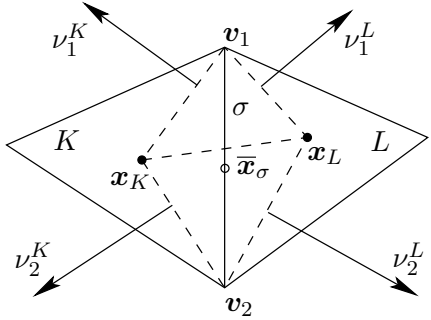
\caption{\label{fig:lmp1}Construction of a monotone scheme on triangular
meshes. The vectors $\nu_i^{K/L}$ have the length of the
segment to which they are orthogonal.}
\end{center}
\end{figure}

A first idea\cite{LEP05} to achieve \eqref{nonlin-2pt} via \eqref{combconv}
on 2D triangular meshes is to compute, for each interior edge $\sigma$
and each $i=1,2$, a constant gradient $\nabla_i u$ on
the triangle $T_i=\V_i \xcv_K\xcv_L$ (see notations in Fig. \ref{fig:lmp1})
by using unknown values $(u_{\V_i},u_K,u_L)$ at this triangle vertices.
These gradients are given by \eqref{eq-grad}
with $\xcv_K$ replaced by $\V_i$ and $\bxs,\bx_{\sigma'}$ replaced
by $\xcv_K,\xcv_L$ and, assuming $\Lambda={\rm Id}$,
the linear conservative fluxes  $F^i_{K,\sigma}$ ($i=1,2$) are then\cite{LEP05,LIP07}
\begin{equation}\label{choixFi}
F^i_{K,\sigma}:=-|\sigma| \nabla_i u\cdot\n_{K,\sigma}=
\frac{|\sigma|}{2|T_i|}\left(u_K \nu_i^L + u_L\nu_i^K - u_{\V_i}(\nu_i^K +\nu_i^L)\right)\cdot
\n_{K,\sigma}
\end{equation}
where $|T_i|$ is the area of triangle $T_i$.
The convex combination \eqref{combconv} is then designed to
eliminate, in $F_{K,\sigma}$,
the term
\[
-\frac{|\sigma|}{2}\left(
\frac{\mu^1_\sigma(u)(\nu_1^K+\nu_1^L)\cdot\n_{K,\sigma}}{|T_1|}u_{\V_1}+
\frac{\mu^2_\sigma(u)(\nu_2^K+\nu_2^L)\cdot\n_{K,\sigma}}{|T_2|}u_{\V_2}\right)
\]
involving $u_{\V_1},u_{\V_2}$ and which prevents this flux from having
the ``2-pt structure'' \eqref{nonlin-2pt}.
As $\nu_1^K+\nu_1^L+\nu_2^K+\nu_2^L=0$, valid choices of the coefficients
are
\begin{equation}\label{choixmu}
\mu^1_\sigma(u)=\frac{u_{\V_2}/|T_2|}{u_{\V_1}/|T_1|+u_{\V_2}/|T_2|}
\mbox{ and }
\mu^2_\sigma(u)=\frac{u_{\V_1}/|T_1|}{u_{\V_1}/|T_1|+u_{\V_2}/|T_2|},
\end{equation}
provided that $u_{\V_1}$ and $u_{\V_2}$ are
both non-negative and not simultaneously equal to $0$ (in this
last case, we can still take $\mu^1_\sigma=\mu^2_\sigma=\frac{1}{2}$). 
Computing these vertex values by convex combinations of the
cell unknowns ensures that they are non-negative whenever
all cell unknowns are non-negative. Two combinations are
suggested in Ref. \refcite{LIP07}, but none of them takes
into account the possible non-smoothness of $\bu$
around discontinuities of $\Lambda$ and the resulting schemes therefore
suffer from a loss of consistency around these discontinuities
(see Remark \ref{rem:convcomb}).

With the choices \eqref{choixFi}-\eqref{choixmu}, it can be
proved that, \emph{provided that $(\xcv_K)_{K\in\mathcal M}$
are at the intersections of the bisectors of the triangles $K\in\mathcal M$}
(this is where the restriction on the
mesh, i.e. that it is made of triangles, comes into play), 
$F_{K,\sigma}$ given by \eqref{combconv} indeed has the ``2pt structure'' \eqref{nonlin-2pt}
with positive coefficients.

\begin{remark}
This construction of fluxes only makes sense if all $u_K$
are non-negative, and the scheme's matrix
$A(U)$ in \eqref{nl-recast} is therefore well defined only for non-negative
cell unknowns. This is not a practical issue as the non-linear
system \eqref{nl-recast} is often solved by iterating an
algorithm of the form $A(U^n)(u^{n+1}_K)_{K\in\mathcal M}=
(B_K(U^n))_{K\in\mathcal M}$ with all components of $B_K(U^n)$
non-negative if all components of $U^n$ are non-negative. By the properties of $A(U)$,
all $u^n_K$ found in these iterations are non-negative and
$A(U^n)$ is therefore well defined.
\end{remark}

The modification of this method\cite{LIP07} for heterogeneous anisotropic tensors $\Lambda$
consists in taking $\xcv_K$ at the intersection of
the bisectors for the $\Lambda_K$-metric of triangle $K$
and in introducing an additional unknown $u_\sigma$ at the
edge midpoint $\bxs$. Four fluxes $F^{i,M}_{K,\sigma}$ are then computed using
gradients in the triangles $\V_i \xcv_M \bxs$
($i=1,2$, $M=K,L$) and the flux continuities $F^{i,K}_{K,\sigma}=
F^{i,L}_{K,\sigma}$ are written to eliminate the unknown
$u_\sigma$ and to obtain two fluxes $F^i_{K,\sigma}$,
which are then used in \eqref{combconv}. New coefficients $\mu^i_\sigma(u)$
are found which eliminate the $u_{\V_i}$ terms
and, thanks to the initial choice of $\xcv_K$,
$F_{K,\sigma}$ has the structure \eqref{nonlin-2pt}.

This method has been extended to 3D tetrahedral meshes in Ref. \refcite{KAP07}
(using convex combinations of three linear fluxes instead of two)
and to general 2D polygonal meshes in Ref. \refcite{LIP07}, albeit in this last case
at the expense of a loss of consistency of the method, especially for
strong anisotropic tensors.

These non-linear 2pt-fluxes methods are not coercive in general
and no convergence proof is provided in the literature. However,
numerical tests show for smooth data a generic order of convergence 2 for the
solution and 1 for its gradient. Some numerical simulations\cite{LIP07}
also confirm that the solution does not satisfy the full discrete
minimum-maximum principle \eqref{disc-minmaxpple} in general:
the approximate solution for $f=0$ may present values beyond
the maximum of the boundary values, and even internal oscillations.

\subsubsection{Polygonal meshes}\label{sec:2ptpoly}

\begin{figure}[h!]
\begin{center}
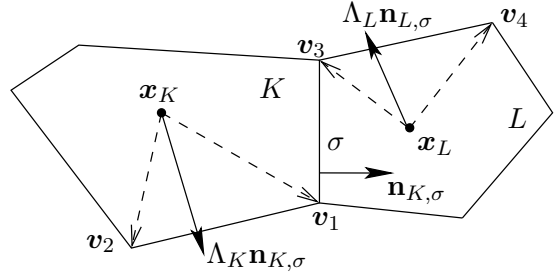
\caption{\label{fig:lmp2}Construction of a monotone scheme on polygonal
meshes.}
\end{center}
\end{figure}

Ref. \refcite{YUA08} presents the construction of consistent 2pt-fluxes \eqref{nonlin-2pt}
on polygonal meshes using similar ideas to Ref. \refcite{LEP05,LIP07}.
The starting point is, for $\sigma\in\mathcal E_K$, to select two vertices $\V_1,\V_2$ of $K$ such that
$\Lambda_K\n_{K,\sigma}$ is in the positive cone generated by $\vect{\xcv_K\V_1}$
and $\vect{\xcv_K\V_2}$ (cf Fig. \ref{fig:lmp2}).

The flux through $\sigma$ outside $K$ can then be approximated by
a positive combination of $-\nabla\bu\cdot \vect{\xcv_K \V_i} \approx \dist(\xcv_K,\V_i)
(u_K-u_{\V_i})$ ($i=1,2$)
and this gives a first numerical flux
$F^1_{K,\sigma}=a_{K,\sigma}(u_K-u_{\V_1})+b_{K,\sigma}(u_K-u_{\V_2})$,
with non-negative coefficients $a_{K,\sigma}$ and $b_{K,\sigma}$. The same
construction from cell $L$ gives a numerical flux outside $K$ (i.e. inside
$L$) 
$F^2_{K,\sigma}=-a_{L,\sigma}(u_L-u_{\V_3})-b_{L,\sigma}(u_L-u_{\V_4})$
with $\V_3,\V_4$ vertices of $L$ and $a_{L,\sigma},b_{L,\sigma}$ non-negative.
The total flux is then obtained as in Refs. \refcite{LEP05,LIP07} by
a convex combination \eqref{combconv} designed to eliminate the coefficients
of $u_{\V_i}$ and to provide the conservativity of the global flux:
\begin{eqnarray*}
\mu^1_{K,\sigma}(u)&=&\frac{a_{L,\sigma}u_{\V_3}+b_{L,\sigma}u_{\V_4}}
{a_{K,\sigma}u_{\V_1}+b_{K,\sigma}u_{\V_2}+a_{L,\sigma}u_{\V_3}+b_{L,\sigma}u_{\V_4}}
\,,\\
\mu^2_{K,\sigma}(u)&=&
\frac{a_{K,\sigma}u_{\V_1}+b_{K,\sigma}u_{\V_2}}
{a_{K,\sigma}u_{\V_1}+b_{K,\sigma}u_{\V_2}+a_{L,\sigma}u_{\V_3}+b_{L,\sigma}u_{\V_4}}.
\end{eqnarray*}
The resulting flux \eqref{combconv} is well defined provided that all $u_{\V_i}$ are
non-negative (if they are all equal to $0$, we take $\mu^i_{K,\sigma}(u)=1/2$) and
has the ``2pt-structure'' \eqref{nonlin-2pt}.
The vertex values $u_{\V_i}$ are computed using convex combinations of
cell values as in Ref. \refcite{LIP07} or, in case of discontinuity of
$\Lambda$, by writing down the flux conservativity and the continuity of tangential
gradients at the vertices. This last method however sometimes fails to provide
non-negative vertex values $u_{\V_i}$ from non-negative cell values, in which
case a simple convex combination must be used.

As for the methods constructed in Refs. \refcite{LEP05,LIP07}, no proof of
convergence is provided in Ref. \refcite{YUA08} but numerical experiments
shows convergence, with rates 2 for $\bu$ and 1 for the fluxes
for smooth data. However, for strongly anisotropic $\Lambda$, the rate of
convergence for $\bu$ seems reduced, at least at available mesh sizes.

This method has been applied to advection-diffusion equations\cite{WAN12} (for
a constant $\Lambda$), using the same kind of discretisation of the advection
term as in Ref. \refcite{LIP10}, i.e. a higher order method with slopes limiters.

A variant can be constructed\cite{SHE12} using edge unknowns $u_\sigma$ (instead
of vertices unknowns) and eliminating them as in the MPFA O-method.
This process may however produce negative $u_\sigma$'s from non-negative $u_K$'s
and, when this happens, $u_\sigma$ must be computed using a simple
convex combination of $u_K$'s.
Although the number of iterations required to compute the solution
are reduced in Ref. \refcite{SHE12} with respect to Ref. \refcite{YUA08},
it seems much higher than for the methods in Refs. \refcite{DRO11,LEP09}
(see Sec. \ref{sec:minmax}), for which the number of
iterations appears to remain bounded independently
on the mesh size.

The ideas of Ref. \refcite{YUA08} have also been used in Ref. \refcite{LIP09-II,LIP10},
but by expressing $\Lambda_K\n_{K,\sigma}$ as a positive combination of
$\vect{\xcv_K \xcv_{L_i}}$, for some cell or edges $L_1,L_2$, instead
of $\vect{\xcv_K \V_i}$ for some vertices $\V_1,\V_2$. This choice
does not require to interpolate new vertex or edge unknowns,
which is an advantage since such interpolations may lead to inaccuracies if not well chosen\cite{LIP07}.
However, when $\Lambda$ is discontinuous
across an edge, the cell centres on each side must be moved according to the heterogeneity
of $\Lambda$ (in such a way that \eqref{cond-orth2} holds for this edge).
As a consequence, the method is applicable only if each cell has
at most one edge across which $\Lambda$ is discontinuous, which
restricts the number and positions of diffusion jumps.
This method has been extended to general 3D polyhedral meshes in 
Refs. \refcite{DAN09,NIK10}.

\subsection{Non-linear multi-point fluxes: MMP schemes}
\label{sec:minmax}

As mentioned above, methods based on the form \eqref{nonlin-2pt} are monotone
but do not satisfy the discrete minimum-maximum principle, mostly because
they do not ensure that $\sum_L \alpha_{K,L}(U)\ge \sum_L \beta_{K,L}(U)$.
It is however possible to construct, on generic 3D meshes, non-linear
MMP schemes provided the fluxes are computed
using a multi-point formula. More precisely, if
\begin{equation}\label{nonlin-mpfluxes}
F_{K,\sigma}=\sum_{Z\in V(K)}\tau_{K,Z}(U)(u_K-u_Z)
\end{equation}
with $V(K)$ a set of cells or edges and $\tau_{K,Z}(U)\ge 0$ ($>0$ whenever
$Z$ is a cell or edge around $K$), then a straightforward
adaptation of the proof in Remark \ref{rem:monvf2} shows that
the resulting scheme satisfies the discrete minimum-maximum principle
\eqref{disc-minmaxpple} (this proof, as mentioned in Remark \ref{rem:monvf2},
demonstrates in fact that the scheme is non-oscillating).
The key element is that \eqref{nonlin-mpfluxes}
ensures that, whenever all cell values are equal, the fluxes
are equal to $0$ or have a sign opposite to the sign of $\udir$
(this is not certain with \eqref{nonlin-2pt}).

A first scheme in this direction is proposed in Ref. \refcite{BER05},
for isotropic diffusion and under restrictive assumptions on the mesh
(made of simplices), such that there exists cell points $(\xcv_K)_{K\in\mathcal M}$
satisfying the orthogonality condition \eqref{cond-orth2}.
For such equations and meshes, the TPFA method can be applied
but the interest of the method in Ref. \refcite{BER05} resides in the
fact that it produces order 2 approximations \emph{of the cell averages of
$\bu$} (the TPFA method would produce order 2 approximations of $(\bu(\xcv_K))_{K\in\mathcal M}$,
where $(\xcv_K)_{K\in\mathcal M}$ are not at cell barycentres).
Nonetheless, the particular convex combinations ideas of Ref. \refcite{BER05}
have been used to construct MMP schemes on triangular
meshes\cite{LEP08,LEP09}, construction then generalised to 
generic 2D or 3D meshes in Ref. \refcite{DRO11}.

\begin{figure}[h!]
\begin{center}
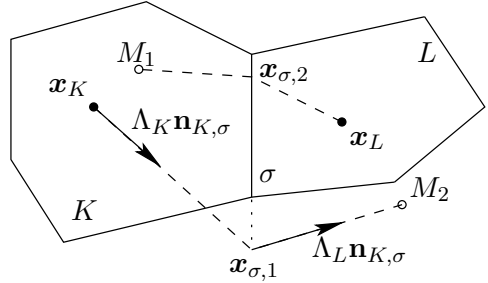
\caption{\label{fig:lmp3}Construction of an MMP scheme.
$\xcv_{\sigma,1}$ is at the intersection of $\xcv_K+[0,\infty)\Lambda_K\n_{K,\sigma}$
and of the line/plane containing $\sigma$. $M_2$ is on the half-line
$\xcv_{\sigma,1}+[0,\infty)\Lambda_L\n_{K,\sigma}$.}
\end{center}
\end{figure}

With the notations in Fig. \ref{fig:lmp3}, the scheme in Ref. \refcite{DRO11}
starts from the two consistent fluxes outside $K$:
\[
\widetilde{F}^1_{K,\sigma}=|\Lambda_K\n_{K,\sigma}|\frac{u_K-u_{\sigma,1}}{\d(\xcv_K,\xcv_{\sigma,1})}\,,\quad \widetilde{F}^2_{K,\sigma}=|\Lambda_L\n_{K,\sigma}|\frac{u_{\sigma,1}-u_{M_2}}{\d(\xcv_{\sigma,1},M_2)}
\]
where $u_{M_2}$ and $u_{\sigma,1}$ are values at $M_2$ and $x_{\sigma,1}$ respectively.
Writing the conservativity
of these fluxes allows us to eliminate $u_{\sigma,1}$ and to get a 
linear conservative flux $F^1_{K,\sigma}=a^1_{K,\sigma}(u_K-u_{M_2})$ with $a^1_{K,\sigma}\ge 0$.
Expressing $u_{M_2}$ as a convex combination of cell unknowns, in such a way
that $u_L$ appears with a non-zero coefficient in this combination (this is
always possible), we then get
\begin{equation}\label{mmp1}
F^1_{K,\sigma}=\alpha^1_{K,\sigma}(u_K-u_L)+G^1_{K,\sigma}\mbox{ with }G^1_{K,\sigma}=
\sum_{M}\beta^1_{K,M}(u_K-u_M)
\end{equation}
with $\alpha^1_{K,\sigma}>0$ and $\beta_{K,M}\ge 0$.
The same construction from cell $L$ gives a flux outside $K$ (i.e. inside $L$)
\begin{equation}\label{mmp2}
F^2_{K,\sigma}=\alpha^2_{L,\sigma}(u_K-u_L)+G^2_{L,\sigma}\mbox{ with }G^2_{L,\sigma}=
\sum_{M}\beta^2_{L,M}(u_M-u_L).
\end{equation}
Following Ref. \refcite{BER05},
a convex combination \eqref{combconv} of these two fluxes is then chosen in order to
eliminate the ``bad'' terms with respect to \eqref{nonlin-mpfluxes}, i.e.
$G^2_{L,\sigma}$:
\begin{equation}\label{choixcombconv}
\mu^1_{K,\sigma}(u)=\frac{|G^2_{L,\sigma}|}{|G^1_{K,\sigma}|+|G^2_{L,\sigma}|}
\,,\quad \mu^2_{K,\sigma}(u)=\frac{|G^1_{K,\sigma}|}{|G^1_{K,\sigma}|+|G^2_{L,\sigma}|}
\end{equation}
(once again, these coefficients are chosen equal to $1/2$ if their denominator
vanishes). By studying separate cases depending on the sign of $G^1_{K,\sigma}G^2_{L,\sigma}$,
we can then see that
$F_{K,\sigma}$ defined by \eqref{combconv}, \eqref{mmp1},
\eqref{mmp2} and \eqref{choixcombconv} always satisfies \eqref{nonlin-mpfluxes},
whatever the values (positive or negative) of the cell unknowns.

\begin{remark}
More freedom is possible on the decompositions in \eqref{mmp1} and \eqref{mmp2},
provided that the global non-linear flux $F_{K,\sigma}$ is
continuous with respect to $u$.
This is ensured\cite{DRO11} if we take $\alpha^1_{K,\sigma}=\alpha^2_{L,\sigma}$
(always possible, upon moving part of the term $u_K-u_L$
in \eqref{mmp1} and \eqref{mmp2} into $G^1_{K,\sigma}$ and $G^2_{L,\sigma}$).
\end{remark}

This method is not necessarily coercive. However, under some coercivity
assumptions (which seem satisfied in numerical tests), a rigorous proof of
convergence is given in Ref. \refcite{DRO11} without regularity assumptions
on the data, drawing on the fact that
the global flux is a convex combination of linear fluxes and adapting the
analysis technique developed in Ref. \refcite{AGE10-2}. This is,
to our best knowledge, the first proof of convergence of an MMP scheme.
Numerical results show a general order 2 convergence for $u$
and, of course, the absence of spurious oscillations in the solution.

\begin{remark}[Choice of convex combination for $u_{M_i}$]\label{rem:convcomb}
In case of jumps of $\Lambda$, numerical tests\cite{DRO11} show that if 
$u_{M_i}$ is computed from cell unknowns on both sides of a discontinuity of $\Lambda$ then the order
of the scheme can be reduced (and the number of Picard iterations
to compute the approximate solution increases significantly).
In many applications, it is however always possible to choose $M_i$ such that
$u_{M_i}$ can be computed using cell unknowns all in a same zone of smoothness
of $\Lambda$.
\end{remark}

The ideas developed for ``2pt non-linear fluxes'' (see Section \ref{sec:nl2pt})
have also been combined with the convex combination
\eqref{choixcombconv} used in Refs. \refcite{BER05,DRO11,LEP05} 
to produce minimum-maximum preserving schemes on 2D polygonal meshes.
In Ref. \refcite{SHE11}, the ideas of Ref. \refcite{YUA08} (replacing vertex unknowns by
edge unknowns) are used to built an MMP method, in which
edge unknowns are interpolated from cell unknowns
by writing a particular flux conservativity which takes into account
the possible jumps of $\Lambda$.

Under an assumption which slightly limits the mesh's skewness and the tensor's
anisotropy, Ref. \refcite{LIP12} draws on the core idea of Ref. \refcite{LIP09-II}
(expressing $\Lambda_K\n_{K,\sigma}$ as a positive combination of $\vect{\xcv_K\xcv_{L_i}}$ for some
cell or edges $L_i$) to produce an MMP scheme
on 2D polygonal meshes. Using cell unknowns rather than interpolating new vertex
or edge unknowns ensures that the stencil of the linear systems
solved at each Picard iteration is as small as the stencil of the TPFA
method (with the trade-of that the fluxes are only conservative
at the limit of these non-linear iterations).
Contrary to Ref. \refcite{LIP09-II}, the method in
Ref. \refcite{LIP12} also does not move cell centres on each side
of an edge across which $\Lambda$ is discontinuous,
but rather makes use of the harmonic interpolation introduced
in Ref. \refcite{AGE09} (see Remark \ref{rem:mix}) to compute the flux through
these edges. The usage of this harmonic interpolation however leads
to a reduced accuracy if the mesh or the tensor are too skewed.

\subsection{MMP schemes by non-linear corrections of linear schemes}

None of the monotone or MMP method presented in the previous sections 
is unconditionally coercive. It turns out that the most efficient
way to construct MMP \emph{and coercive} methods
is not to design a whole new method, but to take existing linear
coercive methods and to devise a non-linear modification of them,
which preserves its coercivity while adding the discrete minimum-maximum
principle.

Let us consider a cell-centred linear scheme \eqref{bf}-\eqref{cf} which is coercive
(it satisfies in particular \eqref{coer-bd}). Assume that, for this scheme,
\[
A_K(u):=\sum_{\sigma\in\mathcal E_K}F_{K,\sigma}=
\sum_{Z\in V(K)}a_{K,Z}(u_K-u_Z)
\]
for some possibly negative $a_{K,Z}$ and $V(K)$ a set of cells or boundary edges
such that, for two cells $(K,Z)$, $Z\in V(K)$ if and only if $K\in V(Z)$.
The scheme is thus written: for all $K\in\mathcal M$, $A_K(u)=\int_K f(x)\d x$.
Then a coercive MMP scheme
can be obtained\cite{LEP10,CAN13} by writing $S_K(u)=\int_K f(x)\d x$ for all $K\in\mathcal M$,
where
\[
\begin{array}{l}
\dsp S_K(u)=A_K(u)+\sum_{Z\in V(K)}\beta_{K,Z}(u)(u_K-u_Z)\\
\dsp \mbox{with }\quad
\beta_{K,Z}(u)\ge \frac{|A_K(u)|}{\sum_{Y\in V(K)}|u_K-u_Y|}
\end{array}
\]
(``$\ge$'' is replaced with ``$>$'' if $Z$ is a neighbouring cell or
edge of $K$, and if $\sum_{Y\in V(K)}|u_K-u_Y|=0$ then we only need $\beta_{K,Z}(u)\ge 0$;
this condition on $\beta_{K,Z}$ is only an example, see Ref. \refcite{CAN13}).
If $\beta_{K,Z}(u)=\beta_{Z,K}(u)$ for any cells $K,Z$,
then the modified scheme is indeed a FV method:
non-linear conservative fluxes $F_{K,\sigma}'(u)$ can be found such that
$S_K(u)=\sum_{\sigma\in\mathcal E_K}F_{K,\sigma}'(u)$.

It is obvious from the symmetry of $\beta_{K,Z}(u)$ that the corrected
scheme retains the coercivity property \eqref{coer-bd} of the original scheme.
It can also be proved that, if the original scheme is consistent in the
sense of FV methods, then the modified scheme converges as the mesh size tends
to $0$, under assumptions on the approximations not formally proved
but holding well in numerical tests. 

These numerical tests show astonishing improvements
of the $L^2$ error when using the non-linear correction
(sometimes\cite{LEP13-p} by a factor 10,000 in case of an anisotropy ratio of $10^6$).
This correction however appears to degrade the order of convergence to
1 and is therefore outperformed by the original order 2 linear
scheme on very thin meshes (sometimes at a size which is nevertheless beyond
computational capacities). The reason for this reduction of
convergence rate is not well understood, but it is worth
mentioning that, even for linear FV schemes, the convergence order 2
on $\bu$ is mostly only \emph{noticed} on numerical tests and
not proved in general.
The consequence is that non-linear corrections should only
be applied for coarse meshes and strongly anisotropic diffusion tensors
for which the original scheme provides physically unacceptable solutions.

This correction technique has been adapted in Ref. \refcite{LEP12} to methods
involving cell and edge unknowns.

\section{Conclusion}\label{sec:concl}

We presented and gave a review of some recent FV methods for diffusion
equations, focusing on the capacity of the methods to
be applicable on generic meshes and to reproduce two important
properties of the continuous equation: coercivity, which
ensures the stability of the scheme and allows one to carry out
convergence proofs under realistic assumptions, and minimum and maximum principles,
which ensure physically acceptable solutions in case of strong anisotropy.

This review is of course partial and much more could be written
on FV methods for \eqref{base}, for example about the comparison
of their respective numerical behaviours -- see e.g. the two comprehensive benchmarks
of Refs. \refcite{EYM12-2,HER08}. Other methods or topics of interest regarding the discretisation
of \eqref{base} are worth mentioning:
\begin{itemize}
\item vertex-centred MPFA O-methods\cite{EDW02,EDW10,EDW11,PAL12},
\item Finite Volume Element methods\cite{CAI91,CAI91-2,EWI02},
based on Finite Element spaces with vertex unknowns
and flux balances on dual meshes around vertices,
\item studies of relationships between FV and Finite Element methods, or mixing
of ideas between different families of methods\cite{VOH06,VOH13,YOU04,WHE06},
\item Gradient Schemes\cite{DRO12,EYM12,EYM11-2,EYM13-2,EYM11},
a generic framework (including HMM methods and some MPFA and DDFV schemes, as well as non-FV
methods) for the convergence analysis of numerical methods on numerous models,
\item the recent review of Ref. \refcite{DIP13-2} on numerical methods in geosciences.
\end{itemize}

The overall conclusion of this review is that currently there is
no miraculous method which provides an excellent
solution in all circumstances. The various numerical
methods available for \eqref{base} should be considered as a kit of clever
techniques which can be adapted and re-used in particular situations.
The ideas behind the methods are as important as the methods 
themselves.

Let us close this study with an open question. For the TPFA scheme,
the flux balance \eqref{bf} can be written
\begin{equation}\label{peculiar}
\sum_{L\in\mathcal M}\tau_{K,L}(u_K-u_L)=\int_K f(x)\d x,
\end{equation}
with $\tau_{K,L}=\tau_{L,K}$ non-negative and such that the method
is coercive.
This structure allows one, by using non-linear functions of the solution
as test functions, to prove \emph{a priori} estimates and analyse the convergence
of the TPFA scheme for non-coercive convection-diffusion equations\cite{DRO02,DRO03-2,CHA11},
hyperbolic-parabolic equations\cite{AND10,EYM02},
equations with Radon measures\cite{GAL99,DRO03} (used to model
wells in reservoirs), or chemotaxis problems\cite{FIL06}.
To date, it is not known how to design a method that can be written
\eqref{peculiar} for any mesh and tensor (as separately noticed in Ref. \refcite{EYM13-3}),
or how to adapt the afore mentioned \emph{a priori} estimate techniques
to schemes not having this structure...

\section*{Acknowledgment}

The author would like to thank the following colleagues, whose comments
helped improve this paper: D. Di Pietro, M.G. Edwards,
R. Eymard, T. Gallou\"et, F. Hermeline, K. Lipnikov, M. Shashkov, D. Svyatskiy
and Yu. Vassilevski.
Special thanks to B. Andreianov, R. Herbin, C. Le Potier and G. Manzini for their thorough reading
and feedback.

\bibliographystyle{abbrv}
\bibliography{m3as-review}

\end{document}

%% file: fig-mesh.pdf_t
\begin{picture}(0,0)%
\includegraphics{fig-mesh.pdf}%
\end{picture}%
\setlength{\unitlength}{4144sp}%
\begingroup\makeatletter\ifx\SetFigFont\undefined%
\gdef\SetFigFont#1#2#3#4#5{%
  \reset@font\fontsize{#1}{#2pt}%
  \fontfamily{#3}\fontseries{#4}\fontshape{#5}%
  \selectfont}%
\fi\endgroup%
\begin{picture}(3144,1475)(472,-254)
\put(1163,610){\makebox(0,0)[lb]{\smash{{\SetFigFont{10}{12.0}{\rmdefault}{\mddefault}{\updefault}{\color[rgb]{0,0,0}$K$}%
}}}}
\put(2655,676){\makebox(0,0)[lb]{\smash{{\SetFigFont{10}{12.0}{\rmdefault}{\mddefault}{\updefault}{\color[rgb]{0,0,0}$L$}%
}}}}
\put(1891,254){\makebox(0,0)[lb]{\smash{{\SetFigFont{10}{12.0}{\rmdefault}{\mddefault}{\updefault}{\color[rgb]{0,0,0}$\sigma$}%
}}}}
\put(2251,434){\makebox(0,0)[lb]{\smash{{\SetFigFont{10}{12.0}{\rmdefault}{\mddefault}{\updefault}{\color[rgb]{0,0,0}$\mathbf{n}_{K,\sigma}$}%
}}}}
\put(676,164){\makebox(0,0)[lb]{\smash{{\SetFigFont{10}{12.0}{\rmdefault}{\mddefault}{\updefault}{\color[rgb]{0,0,0}$\xcv_M$}%
}}}}
\put(901,-106){\makebox(0,0)[lb]{\smash{{\SetFigFont{10}{12.0}{\rmdefault}{\mddefault}{\updefault}{\color[rgb]{0,0,0}$M$}%
}}}}
\put(2601,352){\makebox(0,0)[lb]{\smash{{\SetFigFont{10}{12.0}{\rmdefault}{\mddefault}{\updefault}{\color[rgb]{0,0,0}$\xcv_L$}%
}}}}
\put(1707,721){\makebox(0,0)[lb]{\smash{{\SetFigFont{10}{12.0}{\rmdefault}{\mddefault}{\updefault}{\color[rgb]{0,0,0}$\xcv_K$}%
}}}}
\put(3601,704){\makebox(0,0)[lb]{\smash{{\SetFigFont{10}{12.0}{\rmdefault}{\mddefault}{\updefault}{\color[rgb]{0,0,0}$\Omega$}%
}}}}
\end{picture}%

%% file: fig-mpfa.pdf_t
\begin{picture}(0,0)%
\includegraphics{fig-mpfa.pdf}%
\end{picture}%
\setlength{\unitlength}{4144sp}%
\begingroup\makeatletter\ifx\SetFigFont\undefined%
\gdef\SetFigFont#1#2#3#4#5{%
  \reset@font\fontsize{#1}{#2pt}%
  \fontfamily{#3}\fontseries{#4}\fontshape{#5}%
  \selectfont}%
\fi\endgroup%
\begin{picture}(3526,2004)(1032,-1288)
\put(2508,-782){\makebox(0,0)[lb]{\smash{{\SetFigFont{12}{14.4}{\rmdefault}{\mddefault}{\updefault}{\color[rgb]{0,0,0}$\xcv_K$}%
}}}}
\put(2450,-337){\makebox(0,0)[lb]{\smash{{\SetFigFont{10}{12.0}{\rmdefault}{\mddefault}{\updefault}{\color[rgb]{0,0,0}$K_\V$}%
}}}}
\put(2746,524){\makebox(0,0)[lb]{\smash{{\SetFigFont{12}{14.4}{\rmdefault}{\mddefault}{\updefault}{\color[rgb]{0,0,0}$\xcv_N$}%
}}}}
\put(3286,389){\makebox(0,0)[lb]{\smash{{\SetFigFont{12}{14.4}{\rmdefault}{\mddefault}{\updefault}{\color[rgb]{0,0,0}$\bx_{\sigma''}$}%
}}}}
\put(2476,119){\makebox(0,0)[lb]{\smash{{\SetFigFont{12}{14.4}{\rmdefault}{\mddefault}{\updefault}{\color[rgb]{0,0,0}$\V$}%
}}}}
\put(2650,-1072){\makebox(0,0)[lb]{\smash{{\SetFigFont{12}{14.4}{\rmdefault}{\mddefault}{\updefault}{\color[rgb]{0,0,0}$K$}%
}}}}
\put(1171, 74){\makebox(0,0)[lb]{\smash{{\SetFigFont{12}{14.4}{\rmdefault}{\mddefault}{\updefault}{\color[rgb]{0,0,0}$\xcv_M$}%
}}}}
\put(3871,-691){\makebox(0,0)[lb]{\smash{{\SetFigFont{12}{14.4}{\rmdefault}{\mddefault}{\updefault}{\color[rgb]{0,0,0}$\V'$}%
}}}}
\put(3151,-217){\makebox(0,0)[lb]{\smash{{\SetFigFont{12}{14.4}{\rmdefault}{\mddefault}{\updefault}{\color[rgb]{0,0,0}$\bxs$}%
}}}}
\put(1947,-216){\makebox(0,0)[lb]{\smash{{\SetFigFont{12}{14.4}{\rmdefault}{\mddefault}{\updefault}{\color[rgb]{0,0,0}$\bx_{\sigma'}$}%
}}}}
\put(3601,119){\makebox(0,0)[lb]{\smash{{\SetFigFont{12}{14.4}{\rmdefault}{\mddefault}{\updefault}{\color[rgb]{0,0,0}$\xcv_L$}%
}}}}
\put(3511,-376){\makebox(0,0)[lb]{\smash{{\SetFigFont{12}{14.4}{\rmdefault}{\mddefault}{\updefault}{\color[rgb]{0,0,0}$\sigma$}%
}}}}
\put(4141,-286){\makebox(0,0)[lb]{\smash{{\SetFigFont{12}{14.4}{\rmdefault}{\mddefault}{\updefault}{\color[rgb]{0,0,0}$L$}%
}}}}
\put(1711,-601){\makebox(0,0)[lb]{\smash{{\SetFigFont{12}{14.4}{\rmdefault}{\mddefault}{\updefault}{\color[rgb]{0,0,0}$\sigma'$}%
}}}}
\put(4051, 74){\makebox(0,0)[lb]{\smash{{\SetFigFont{12}{14.4}{\rmdefault}{\mddefault}{\updefault}{\color[rgb]{0,0,0}$\n_{K,\sigma}$}%
}}}}
\put(2116,-826){\makebox(0,0)[lb]{\smash{{\SetFigFont{12}{14.4}{\rmdefault}{\mddefault}{\updefault}{\color[rgb]{0,0,0}$\nu_{K,\sigma'}$}%
}}}}
\put(3422,-1102){\makebox(0,0)[lb]{\smash{{\SetFigFont{12}{14.4}{\rmdefault}{\mddefault}{\updefault}{\color[rgb]{0,0,0}$\nu_{K,\sigma}$}%
}}}}
\end{picture}%

%% file: fig-ddfv.pdf_t
\begin{picture}(0,0)%
\includegraphics{fig-ddfv.pdf}%
\end{picture}%
\setlength{\unitlength}{4144sp}%
\begingroup\makeatletter\ifx\SetFigFont\undefined%
\gdef\SetFigFont#1#2#3#4#5{%
  \reset@font\fontsize{#1}{#2pt}%
  \fontfamily{#3}\fontseries{#4}\fontshape{#5}%
  \selectfont}%
\fi\endgroup%
\begin{picture}(4167,2191)(706,-4220)
\put(3871,-3031){\makebox(0,0)[lb]{\smash{{\SetFigFont{12}{14.4}{\rmdefault}{\mddefault}{\updefault}{\color[rgb]{0,0,0}$K$}%
}}}}
\put(3691,-3706){\makebox(0,0)[lb]{\smash{{\SetFigFont{12}{14.4}{\rmdefault}{\mddefault}{\updefault}{\color[rgb]{0,0,0}$\xcv_K$}%
}}}}
\put(1396,-3121){\makebox(0,0)[lb]{\smash{{\SetFigFont{12}{14.4}{\rmdefault}{\mddefault}{\updefault}{\color[rgb]{0,0,0}$L$}%
}}}}
\put(1801,-3661){\makebox(0,0)[lb]{\smash{{\SetFigFont{12}{14.4}{\rmdefault}{\mddefault}{\updefault}{\color[rgb]{0,0,0}$\xcv_L$}%
}}}}
\put(721,-3481){\makebox(0,0)[lb]{\smash{{\SetFigFont{12}{14.4}{\rmdefault}{\mddefault}{\updefault}{\color[rgb]{0,0,0}$M$}%
}}}}
\put(1189,-3431){\makebox(0,0)[lb]{\smash{{\SetFigFont{12}{14.4}{\rmdefault}{\mddefault}{\updefault}{\color[rgb]{0,0,0}$\sigma'$}%
}}}}
\put(1036,-3931){\makebox(0,0)[lb]{\smash{{\SetFigFont{12}{14.4}{\rmdefault}{\mddefault}{\updefault}{\color[rgb]{0,0,0}$\xcv_M$}%
}}}}
\put(2476,-4156){\makebox(0,0)[lb]{\smash{{\SetFigFont{12}{14.4}{\rmdefault}{\mddefault}{\updefault}{\color[rgb]{0,0,0}$\V'$}%
}}}}
\put(2611,-2581){\makebox(0,0)[lb]{\smash{{\SetFigFont{12}{14.4}{\rmdefault}{\mddefault}{\updefault}{\color[rgb]{0,0,0}$\V$}%
}}}}
\put(2881,-3211){\makebox(0,0)[lb]{\smash{{\SetFigFont{12}{14.4}{\rmdefault}{\mddefault}{\updefault}{\color[rgb]{0,0,0}$D$}%
}}}}
\put(1891,-2941){\makebox(0,0)[lb]{\smash{{\SetFigFont{12}{14.4}{\rmdefault}{\mddefault}{\updefault}{\color[rgb]{0,0,0}$\mathbf{n}_{K,\sigma}$}%
}}}}
\put(3151,-2536){\makebox(0,0)[lb]{\smash{{\SetFigFont{12}{14.4}{\rmdefault}{\mddefault}{\updefault}{\color[rgb]{0,0,0}$P_\V$}%
}}}}
\put(3196,-3886){\makebox(0,0)[lb]{\smash{{\SetFigFont{12}{14.4}{\rmdefault}{\mddefault}{\updefault}{\color[rgb]{0,0,0}$\mathbf{n}_{\V,\tau}$}%
}}}}
\put(2372,-3405){\makebox(0,0)[lb]{\smash{{\SetFigFont{12}{14.4}{\rmdefault}{\mddefault}{\updefault}{\color[rgb]{0,0,0}$x_\sigma$}%
}}}}
\put(2669,-3734){\makebox(0,0)[lb]{\smash{{\SetFigFont{10}{12.0}{\rmdefault}{\mddefault}{\updefault}{\color[rgb]{0,0,0}$\widehat{\sigma\tau}$}%
}}}}
\put(2947,-3486){\makebox(0,0)[lb]{\smash{{\SetFigFont{12}{14.4}{\rmdefault}{\mddefault}{\updefault}{\color[rgb]{0,0,0}$\tau$}%
}}}}
\put(2600,-3209){\makebox(0,0)[lb]{\smash{{\SetFigFont{12}{14.4}{\rmdefault}{\mddefault}{\updefault}{\color[rgb]{0,0,0}$\sigma$}%
}}}}
\end{picture}%

%% file: fig-lmp1.pdf_t
\begin{picture}(0,0)%
\includegraphics{fig-lmp1.pdf}%
\end{picture}%
\setlength{\unitlength}{4144sp}%
\begingroup\makeatletter\ifx\SetFigFont\undefined%
\gdef\SetFigFont#1#2#3#4#5{%
  \reset@font\fontsize{#1}{#2pt}%
  \fontfamily{#3}\fontseries{#4}\fontshape{#5}%
  \selectfont}%
\fi\endgroup%
\begin{picture}(2544,1936)(1339,-4186)
\put(2701,-2941){\makebox(0,0)[lb]{\smash{{\SetFigFont{12}{14.4}{\rmdefault}{\mddefault}{\updefault}{\color[rgb]{0,0,0}$\sigma$}%
}}}}
\put(3151,-3211){\makebox(0,0)[lb]{\smash{{\SetFigFont{12}{14.4}{\rmdefault}{\mddefault}{\updefault}{\color[rgb]{0,0,0}$\xcv_L$}%
}}}}
\put(1631,-3152){\makebox(0,0)[lb]{\smash{{\SetFigFont{12}{14.4}{\rmdefault}{\mddefault}{\updefault}{\color[rgb]{0,0,0}$K$}%
}}}}
\put(2729,-3296){\makebox(0,0)[lb]{\smash{{\SetFigFont{12}{14.4}{\rmdefault}{\mddefault}{\updefault}{\color[rgb]{0,0,0}$\bxs$}%
}}}}
\put(2656,-4117){\makebox(0,0)[lb]{\smash{{\SetFigFont{12}{14.4}{\rmdefault}{\mddefault}{\updefault}{\color[rgb]{0,0,0}$\V_2$}%
}}}}
\put(2645,-2502){\makebox(0,0)[lb]{\smash{{\SetFigFont{12}{14.4}{\rmdefault}{\mddefault}{\updefault}{\color[rgb]{0,0,0}$\V_1$}%
}}}}
\put(3529,-3158){\makebox(0,0)[lb]{\smash{{\SetFigFont{12}{14.4}{\rmdefault}{\mddefault}{\updefault}{\color[rgb]{0,0,0}$L$}%
}}}}
\put(1970,-3359){\makebox(0,0)[lb]{\smash{{\SetFigFont{12}{14.4}{\rmdefault}{\mddefault}{\updefault}{\color[rgb]{0,0,0}$\xcv_K$}%
}}}}
\put(3660,-3854){\makebox(0,0)[lb]{\smash{{\SetFigFont{12}{14.4}{\rmdefault}{\mddefault}{\updefault}{\color[rgb]{0,0,0}$\nu_2^L$}%
}}}}
\put(1621,-2542){\makebox(0,0)[lb]{\smash{{\SetFigFont{12}{14.4}{\rmdefault}{\mddefault}{\updefault}{\color[rgb]{0,0,0}$\nu_1^K$}%
}}}}
\put(3368,-2539){\makebox(0,0)[lb]{\smash{{\SetFigFont{12}{14.4}{\rmdefault}{\mddefault}{\updefault}{\color[rgb]{0,0,0}$\nu_1^L$}%
}}}}
\put(1398,-3857){\makebox(0,0)[lb]{\smash{{\SetFigFont{12}{14.4}{\rmdefault}{\mddefault}{\updefault}{\color[rgb]{0,0,0}$\nu_2^K$}%
}}}}
\end{picture}%

%% file: fig-lmp2.pdf_t
\begin{picture}(0,0)%
\includegraphics{fig-lmp2.pdf}%
\end{picture}%
\setlength{\unitlength}{4144sp}%
\begingroup\makeatletter\ifx\SetFigFont\undefined%
\gdef\SetFigFont#1#2#3#4#5{%
  \reset@font\fontsize{#1}{#2pt}%
  \fontfamily{#3}\fontseries{#4}\fontshape{#5}%
  \selectfont}%
\fi\endgroup%
\begin{picture}(3264,1684)(1879,-4229)
\put(2656,-3166){\makebox(0,0)[lb]{\smash{{\SetFigFont{12}{14.4}{\rmdefault}{\mddefault}{\updefault}{\color[rgb]{0,0,0}$\xcv_K$}%
}}}}
\put(3376,-3166){\makebox(0,0)[lb]{\smash{{\SetFigFont{12}{14.4}{\rmdefault}{\mddefault}{\updefault}{\color[rgb]{0,0,0}$K$}%
}}}}
\put(2341,-4066){\makebox(0,0)[lb]{\smash{{\SetFigFont{12}{14.4}{\rmdefault}{\mddefault}{\updefault}{\color[rgb]{0,0,0}$\V_2$}%
}}}}
\put(3781,-3481){\makebox(0,0)[lb]{\smash{{\SetFigFont{12}{14.4}{\rmdefault}{\mddefault}{\updefault}{\color[rgb]{0,0,0}$\sigma$}%
}}}}
\put(4321,-3481){\makebox(0,0)[lb]{\smash{{\SetFigFont{12}{14.4}{\rmdefault}{\mddefault}{\updefault}{\color[rgb]{0,0,0}$\xcv_L$}%
}}}}
\put(4816,-2716){\makebox(0,0)[lb]{\smash{{\SetFigFont{12}{14.4}{\rmdefault}{\mddefault}{\updefault}{\color[rgb]{0,0,0}$\V_4$}%
}}}}
\put(4861,-3346){\makebox(0,0)[lb]{\smash{{\SetFigFont{12}{14.4}{\rmdefault}{\mddefault}{\updefault}{\color[rgb]{0,0,0}$L$}%
}}}}
\put(3871,-2716){\makebox(0,0)[lb]{\smash{{\SetFigFont{12}{14.4}{\rmdefault}{\mddefault}{\updefault}{\color[rgb]{0,0,0}$\Lambda_L\n_{L,\sigma}$}%
}}}}
\put(4141,-3751){\makebox(0,0)[lb]{\smash{{\SetFigFont{12}{14.4}{\rmdefault}{\mddefault}{\updefault}{\color[rgb]{0,0,0}$\n_{K,\sigma}$}%
}}}}
\put(3061,-4156){\makebox(0,0)[lb]{\smash{{\SetFigFont{12}{14.4}{\rmdefault}{\mddefault}{\updefault}{\color[rgb]{0,0,0}$\Lambda_K\n_{K,\sigma}$}%
}}}}
\put(3691,-3931){\makebox(0,0)[lb]{\smash{{\SetFigFont{12}{14.4}{\rmdefault}{\mddefault}{\updefault}{\color[rgb]{0,0,0}$\V_1$}%
}}}}
\put(3601,-2896){\makebox(0,0)[lb]{\smash{{\SetFigFont{12}{14.4}{\rmdefault}{\mddefault}{\updefault}{\color[rgb]{0,0,0}$\V_3$}%
}}}}
\end{picture}%

%% file: fig-lmp3.pdf_t
\begin{picture}(0,0)%
\includegraphics{fig-lmp3.pdf}%
\end{picture}%
\setlength{\unitlength}{4144sp}%
\begingroup\makeatletter\ifx\SetFigFont\undefined%
\gdef\SetFigFont#1#2#3#4#5{%
  \reset@font\fontsize{#1}{#2pt}%
  \fontfamily{#3}\fontseries{#4}\fontshape{#5}%
  \selectfont}%
\fi\endgroup%
\begin{picture}(2859,1705)(2284,-4319)
\put(4321,-3481){\makebox(0,0)[lb]{\smash{{\SetFigFont{12}{14.4}{\rmdefault}{\mddefault}{\updefault}{\color[rgb]{0,0,0}$\xcv_L$}%
}}}}
\put(2656,-3931){\makebox(0,0)[lb]{\smash{{\SetFigFont{12}{14.4}{\rmdefault}{\mddefault}{\updefault}{\color[rgb]{0,0,0}$K$}%
}}}}
\put(4726,-2986){\makebox(0,0)[lb]{\smash{{\SetFigFont{12}{14.4}{\rmdefault}{\mddefault}{\updefault}{\color[rgb]{0,0,0}$L$}%
}}}}
\put(3781,-3706){\makebox(0,0)[lb]{\smash{{\SetFigFont{12}{14.4}{\rmdefault}{\mddefault}{\updefault}{\color[rgb]{0,0,0}$\sigma$}%
}}}}
\put(3016,-3391){\makebox(0,0)[lb]{\smash{{\SetFigFont{12}{14.4}{\rmdefault}{\mddefault}{\updefault}{\color[rgb]{0,0,0}$\Lambda_K\n_{K,\sigma}$}%
}}}}
\put(2521,-3166){\makebox(0,0)[lb]{\smash{{\SetFigFont{12}{14.4}{\rmdefault}{\mddefault}{\updefault}{\color[rgb]{0,0,0}$\xcv_K$}%
}}}}
\put(3601,-4246){\makebox(0,0)[lb]{\smash{{\SetFigFont{12}{14.4}{\rmdefault}{\mddefault}{\updefault}{\color[rgb]{0,0,0}$\xcv_{\sigma,1}$}%
}}}}
\put(4096,-4156){\makebox(0,0)[lb]{\smash{{\SetFigFont{12}{14.4}{\rmdefault}{\mddefault}{\updefault}{\color[rgb]{0,0,0}$\Lambda_L\n_{K,\sigma}$}%
}}}}
\put(4681,-3796){\makebox(0,0)[lb]{\smash{{\SetFigFont{12}{14.4}{\rmdefault}{\mddefault}{\updefault}{\color[rgb]{0,0,0}$M_2$}%
}}}}
\put(2926,-2986){\makebox(0,0)[lb]{\smash{{\SetFigFont{12}{14.4}{\rmdefault}{\mddefault}{\updefault}{\color[rgb]{0,0,0}$M_1$}%
}}}}
\put(3781,-3076){\makebox(0,0)[lb]{\smash{{\SetFigFont{12}{14.4}{\rmdefault}{\mddefault}{\updefault}{\color[rgb]{0,0,0}$\xcv_{\sigma,2}$}%
}}}}
\end{picture}%